\documentclass[11pt]{amsart}
\usepackage{amsfonts,latexsym,amsmath,amscd,geometry,graphicx}
\usepackage{amssymb}
\usepackage{latexsym}
\usepackage{xcolor}

\newcommand \nc{\newcommand}
\newtheorem{theorem}{Theorem}[section]
\newtheorem{lemma}[theorem]{Lemma}

\nc{\ba}{\begin{array}}\nc{\ea}{\end{array}}
\nc{\be}{\begin{eqnarray}}\nc{\ee}{\end{eqnarray}}
\nc{\beq}{\begin{equation}}\nc{\eeq}{\end{equation}}
\nc{\bex}{\begin{eqnarray*}}\nc{\eex}{\end{eqnarray*}}
\nc{\btm}{\begin{theorem}} \nc{\etm}{\end{theorem}}
\nc{\blm}{\begin{lemma}} \nc{\elm}{\end{lemma}}
\nc{\R}{\mathbb{R}}  \nc{\ld}{\lambda}
\nc{\va}{\varphi}
\nc{\ve}{\varepsilon}
\def\tu{\tilde{u}}

\def\pf{\noindent{\bf Proof.\quad}}

\numberwithin{equation}{section}

\makeatletter
\@namedef{subjclassname@2020}{%
  \textup{2020} Mathematics Subject Classification}
\makeatother

\allowdisplaybreaks[1]

\begin{document}

\title{BV solutions to a hyperbolic system of balance laws with logistic growth}
% on a Keller-Segel type chemotaxis model}

\author{Geng Chen}
\address[G. Chen]{Department of Mathematics, University of Kansas, Lawrence, KS 66045, U.S.A.}
\email{gengchen@ku.edu}

\author{Yanni Zeng}
\address[Y. Zeng]{Department of Mathematics, University of Alabama at Birmingham, Birmingham, AL 35294, U.S.A}
\email{ynzeng@uab.edu}


\begin{abstract} 
We study BV solutions for a $2\times2$ system of hyperbolic balance laws. We show that when initial data have small total variation on $(-\infty,\infty)$ and small amplitude, and decay sufficiently fast to a constant equilibrium state as $|x|\rightarrow\infty$, a Cauchy problem (with generic data) has a unique admissible BV solution defined globally in time. Here the solution is admissible in the sense that its shock waves satisfy the Lax entropy condition. We also study asymptotic behavior of solutions. In particular, we obtain a time decay rate for the total variation of the solution, and a convergence rate of the solution to its time asymptotic solution. Our system is a modification of a Keller-Segel type chemotaxis model. Its flux function possesses new features when comparing to the well-known model of Euler equations with damping. This may help to shed light on how to extend the study to a general system of hyperbolic balance laws in the future.

\end{abstract}

\maketitle

\section{Introduction}
We consider the Cauchy problem of a $2\times2$ system of hyperbolic balance laws,
\beq\label{KS}
\left\{
\begin{split}
v_t+u_x &=0,\\
u_t+(uv)_x&=u(1-u),
\end{split}
\right.
\eeq
where the initial data satisfies
\beq\label{ID}
(v,u)(x,0)=(v_0,u_0)(x),\qquad \lim_{x\rightarrow \pm \infty}(v_0,u_0)(x)=(0,1).
\eeq

The model is the converted form under inverse Hopf-Cole transformation of a
Keller-Segel type chemotaxis model with logistic growth,
logarithmic sensitivity and density-dependent production/consumption rate. Here we give a brief discussion on its background.

The original model, Othmer-Stevens model \cite{OS,LS},  describes the dynamical behavior of chemotactic movement of random walkers that deposit non-diffusive or slow-moving chemical signals to modify the local environment for succeeding passages. The model reads
\begin{equation}\label{OS}
\begin{cases}
s_t=-\mu us-\sigma s,\\
u_t=Du_{xx}-\chi[u(\ln s)_x]_x,
\end{cases}
\end{equation}
where the unknown functions are $s=s(x,t)$ and $u=u(x,t)$ for the concentration of a chemical signal and density of a cellular population, respectively. We have assumed that the chemical signal is non-diffusive. Meanwhile, the system parameters have the following meaning:
\begin{itemize}
\item $\mu\neq0$: coefficient of density-dependent production/consumption rate of chemical signal;
\item $\sigma\ge0$:  natural degradation rate of chemical signal;
\item $D\ge0$: diffusion coefficient of cellular population;
\item $\chi\neq0$: coefficient of chemotactic sensitivity.
\end{itemize}

Mathematical properties of \eqref{OS}, along with its companion with chemical diffusion, have been studied abundantly in recent years. These include, for instance, global well-posedness \cite{FFH,GXZZ}, asymptotic stability of equilibria \cite{1d3,1d1,1d4,CLW}, existence and stability of traveling wave solutions \cite{TW2,TW8,TW1}, and vanishing chemical diffusivity limit \cite{DL2,1d5}. Also see references therein.

We assume
\begin{equation}\label{1.4}
\chi\mu>0,
\end{equation}
which includes two scenarios: $\chi>0$ and $\mu>0$,  or $\chi<0$ and $\mu<0$. The former is interpreted as cells are attracted to and consume the chemical. On the other hand, the latter describes cells depositing the chemical to modify the local environment for succeeding passages \cite{OS}. Mathematically, \eqref{1.4} allows us to convert \eqref{OS} (with $D=0$) to a hyperbolic system  in biologically relevant regimes \cite{ZZ1}. 

Next we append \eqref{OS} by logistic growth of the cellular population  to account for limited resources:
\begin{equation}\label{1.5}
\begin{cases}
s_t=-\mu us-\sigma s,\\
u_t=Du_{xx}-\chi[u(\ln s)_x]_x+au(1-\frac uK),
\end{cases}
\end{equation}
where  $a>0$ is the  natural growth rate of cellular population and $K>0$ is the typical carrying capacity. 

The logarithmic sensitivity function in \eqref{OS} or \eqref{1.5} is based on the assumption that the detection of chemical signal by cellular population follows the Weber-Fechner law. The singularity can be removed via the inverse Hopf-Cole transformation \cite{LS}:
\begin{equation}\label{1.6}
v=(\ln s)_x=\frac{s_x}s.
\end{equation}
Under the new variables $v$ and $u$, the reaction-diffusion-advection system \eqref{1.5} becomes a system of hyperbolic-parabolic balance laws:
\begin{equation}\label{1.7}
\begin{cases}
v_t+\mu u_x=0,\\
u_t+\chi(uv)_x=Du_{xx}+au(1-\frac uK).
\end{cases}
\end{equation}

Under the assumption \eqref{1.4},  \eqref{1.7} can be further simplified by rescaling:
\begin{equation}\label{1.8}
\tilde t=\chi\mu Kt, \qquad \tilde x=\sqrt{\chi\mu K}x,\qquad \tilde v=\mathrm{sign}(\chi)\sqrt{\frac{\chi}{\mu K}}v,
\qquad\tilde u=\frac uK.
\end{equation}
After dropping the tilde accent, we arrive at
\begin{equation}\label{1.9}
\begin{cases}
v_t+u_x=0,\\
u_t+(uv)_x=Du_{xx}+ru(1- u),
\end{cases}
\end{equation}
where
\begin{equation}\label{1.10}
r=\frac{a}{\chi\mu K}>0.
\end{equation}
As our main goal in this paper is to establish global existence of BV solutions for hyperbolic balance laws (to be discussed below), we set $D=0$ and without loss of generality, $r=1$. This gives us \eqref{KS}.

We consider Cauchy problem around a constant equilibrium state. To be an equilibrium state, $u=0$ or $u=1$ while $v$ can be any constant. For stability we take $u=1$.  For physical relevance (assuming the existence of $\lim s$ as ${x\rightarrow\pm\infty}$) and by \eqref{1.6}, we take $v=0$ in the equilibrium state. These give us the setup in initial data \eqref{ID}.

The existence framework on small BV solutions for a general system of hyperbolic balance laws in the form of
\begin{equation}\label{1.11}
\left\{
\begin{split}
\partial_t V+\partial_x P(V,W) =0\\
\partial_t W+\partial_x Q(V,W) +\Omega(V,W)W=0
\end{split}
\right.
\end{equation}
has been established by a sequence of works of Dafermos including \cite{D06,D13,D14,D15,Dafermos}, following the pioneering work of  Dafermos and Hsiao in \cite{DH82}. Here $V\in \mathbb{R}^k$,  $W\in \mathbb{R}^l$. The main idea in \cite{DH82} is to modify the random choice method introduced by Glimm \cite{Glimm} for hyperbolic conservation laws. More precisely, one applies an operator splitting technique in the random choice method to cope with the damping term. 

Unfortunately, in current existence results for the general system,  one still needs to assume that the total mass of $V$ is zero by imposing the condition
\beq\label{v0con}
\int_{-\infty}^\infty V_0(x)dx=0
\eeq
on the initial data, or assume that the equilibrium state is $L^1$-stable, see for instance \cite{D13,D15,D06}. There are important special cases where the restriction \eqref{v0con} can be removed. They are the p-system and its variety \cite{D95,D14}:
\beq\label{psys}
\left\{
\begin{split}
v_t-u_x &=-\alpha u^2,\\
u_t+p(v)_x&=-u.
\end{split}
\right.
\eeq
Here $p$ is a given smooth function with $p(0)=0$ and $p'(0)=-1$, and $\alpha$ is an arbitrarily given constant (including zero). Another variety is 
\beq\label{psys+}
\left\{
\begin{split}
v_t-u_x &=0,\\
u_t+p(v)_x&=-u+f(v),
\end{split}
\right.
\eeq
where $p$ and $f$ satisfy the subcharacterisitc condition
$$
-p'(v)>[f'(v)]^2,
$$
see \cite{D06}.

The goal of this paper is to study small BV solutions to \eqref{KS} without assuming zero mass on $v$.  This includes global existence and long time behavior.
 Our intention is that by studying a different, physically relevant model, one gains new insights into technical difficulties associated with the general system \eqref{1.11} when the restriction \eqref{v0con} is removed.

Although \eqref{KS} looks somewhat similar to the p-system, \eqref{psys} with $\alpha=0$, and \eqref{psys+}, we observe two significant differences. 
The first one is that the nonlinear flux functions in $\eqref{psys}_2$ and $\eqref{psys+}_2$ are $p(v)$, depending only on $v$, while the counterpart in $\eqref{KS}_2$ is $uv$, depending on both $v$ and $u$. The second one is the damping in the second equation. In \eqref{psys} and \eqref{psys+} the damping terms are linear  in $u$ while in \eqref{KS} the logistic growth induces an extra higher order term in $u$, see \eqref{KS2} below.

In a recent study of diffusive contact waves for \eqref{KS}, the first difference indeed alters the components in the time asymptotic ansatz and results in slower decay to the background wave while the second one causes some technical difficulties when handling transitional end-states \cite{Zeng3, Zeng4}.  Similarly, when obtaining key  $L^1$ estimates here,  the new form of flux function makes our analysis of \eqref{KS} much more complicated. Besides, when handling the total variation (another key estimate), results from \cite{DH82} can be applied to \eqref{psys} directly \cite{D14}. By contrast, here we need to combine techniques of operator splitting from both \cite{DH82} and \cite{D06}, due to the dependence of the flux function on both unknown variables.

The current paper is an addition to \cite{D14} in contributing our understanding beyond the p-system \cite{D95} for small but generic BV solutions. In particular, our result may shed light on how to remove the stringent restriction \eqref{v0con} in the BV existence theory for the general system of balance laws \eqref{1.11} in a future research.

Now we give details of our result. For convenience, we set $\tu=u-1$ and $\tu_0=u_0-1$ to move the equilibrium state to $(0,0)$. Thus \eqref{KS}-\eqref{ID} become
\beq\label{KS2}
\left\{
\begin{split}
v_t+\tu_x &=0,\\
\tu_t+\big((\tu+1) v\big)_x&=-\tu(1+\tu),
\end{split}
\right.
\eeq
where the initial data satisfy
\beq\label{ID2}
(v,\tu)(x,0)=(v_0,\tu_0)(x),\qquad \lim_{x\rightarrow \pm \infty}(v_0,\tu_0)(x)=(0,0).
\eeq

We are interested in BV solutions of the initial value problem \eqref{KS}-\eqref{ID}, or equivalently \eqref{KS2}-\eqref{ID2}, with small initial data. More precisely, we set
\beq\label{IDTV}
{TV}v_0(\cdot)+TV\tu_0(\cdot)=\delta.
\eeq
Assuming the initial data approach zero sufficiently fast as $x\rightarrow\pm\infty$, say, $(v_0,\tu_0)(x)=O(|x|^{-r})$ with $r>\frac32$, we also set
\beq\label{IDint}
\int_{-\infty}^\infty (1+x^2)[v_0^2(x)+\tu_0^2(x)]\, dx=\sigma^2.
\eeq
Our basic assumption is that $\delta$ and $\sigma$ are small.

We recall that a BV solution $(v,\tu)$ of \eqref{KS2}-\eqref{ID2} on $(-\infty,\infty)\times [0,T)$ for any $T>0$ with $\max|(v,\tu)|< \rho_0$ for some small positive constant $\rho_0$ is called admissible when  its shock waves satisfy the Lax entropy condition. The following is our main theorem, which gives both global existence of a unique admissible BV solution and its large time behavior.

\begin{theorem}\label{main_thm}
There exist positive constants $\delta_0$ and $\sigma_0$ such that if \eqref{IDTV} and \eqref{IDint} hold with $\delta<\delta_0$ and $\sigma<\sigma_0$, then the Cauchy problem \eqref{KS}-\eqref{ID}, with $(v_0,u_0-1)(x)=O(|x|^{-r})$ as $x\rightarrow\pm\infty$ for some $r>\frac32$,  has a unique admissible BV solution $(v,u)$ defined on $(-\infty,\infty)\times[0,\infty)$. For $0\le t<\infty$, the solution has the following decay estimates,
\beq\label{main_thm_1}
TV v(\cdot,t)+TV u(\cdot,t)\leq b\sigma (t+1)^{-\frac{1}{4}}+b\delta e^{-\nu t}
\eeq
and
\beq\label{main_thm_1a}
\int_{-\infty}^\infty
\Big(|v(x,t)-\theta(x,t)|+|u(x,t)-1|\Big)\, dx\leq c\sigma(t+1)^{-\frac{1}{4}},
\eeq
where $b$, $c$ and $\nu$ are positive constants  independent of the initial data, and 
\beq\label{main_thm_1b}
\theta(x,t)=\frac{M}{\sqrt{4\pi(t+1)}} e^{-\frac{x^2}{4(t+1)}},\qquad
M=\int_{-\infty}^\infty v_0(x)\,dx.
\eeq
\end{theorem}

We note that $\theta(x,t)$ satisfies the heat equation,
\begin{equation}\label{heat}
\theta_t=\theta_{xx}.
\end{equation}

Theorem \ref{main_thm} is parallel to Theorem 1.1 in \cite{D14} for \eqref{psys}. Our approach, however, has several important differences than the one in  \cite{D14}.

For key $L^1$ estimates, while we generally follow the road map in \cite{D14} (but with a different set of technical difficulties associated with \eqref{KS}), we bypass an auxiliary system introduced in \cite{D14}. Precisely,  Chapman-Enskog expansion is used in \cite{D14} to reduce the second equation in \eqref{psys}. This gives rise to an auxiliary system with its solution $(\hat v, \hat u)$. Properties of $(\hat v, \hat u)$ are derived and the approximation of $(v,u)$ by $(\hat v, \hat u)$ is studied. Finally, $\hat v$ is further approximated by a heat kernel similar to our $\theta$ in \eqref{main_thm_1b}.
In this paper we perform Chapman-Enskog expansion to the whole system \eqref{KS}, and construct the time asymptotic solution $\theta$ for $v$. We approximate $v$ by $\theta$ directly.

Our approach has several advantages. Since $\theta$ is explicitly formulated in \eqref{main_thm_1b}, its properties are straightforward. It is not that obvious how the auxiliary system and its solution $(\hat v, \hat u)$ 
in \cite{D14} can be extended to a general system \eqref{1.11}. Bypassing such a step reduces technical difficulties associated with constructing the approximate solution and studying its properties. Besides, Chapman-Enskog expansion has been done to \eqref{1.11}, and a time asymptotic solution has been constructed in a systematic way, using heat kernels and Burgers kernels \cite{ZC}. The explicit formulation of the asymptotic solution gives us the needed properties of it right away.  We hope that our new approach 
moves us one step forward to the establishment of BV existence for the general system \eqref{1.11}.

Another key step is to obtain bounds on total variations to extend a local solution to a global one. After redistributing the damping mechanism, one arrives at a new nonhomogeneous hyperbolic system of balance laws. For \eqref{psys} considered in \cite{D14}, the new flux function has good regularity, and hence results from \cite{DH82} apply.
In our case, due to the fact that the flux function in $\eqref{KS}_2$ (or in $\eqref{KS2}_2$) contains both unknown variables, the new flux function has lower regularity and results from \cite{DH82} do not apply. The regularity is the same as those considered in \cite{D06,D13,D15} but here we have extra inhomogeneity under nonzero mass situation, and hence results in those works do not apply either. Therefore, a new approach that combines techniques of operator splitting from \cite{DH82} and \cite{D06} is needed in our work. In this case, the explicit formulation of the asymptotic solution is crucial in our key process of controlling and estimating total variations. At the same time, we incorporate in the process the time decay rate of the $L^1$-norm of the solution to obtain a time decay rate for the total variation, which is not available in \cite{D06,D13,D15}.

The paper is divided into four sections. In section 2 we give some basic setup of our analysis. Three key lemmas on time decay, including $L^2$ and $L^1$ estimates, will be given in section 3. In section 4, we estimate total variations of solutions  to prove our main theorem.

\section{Basic setup}

Since we are considering admissible BV solutions of \eqref{KS}, \eqref{ID} near $(0,1)$, or equivalently, admissible BV solutions of \eqref{KS2}, \eqref{ID2} near $(0,0)$, throughout this paper our solutions and approximate solutions of \eqref{KS2} take values in a small neighborhood of the origin. 

Form \eqref{KS2}, the Jacobian matrix of the flux function has two distinct, real eigenvalues
\begin{equation}\label{2.1}
\lambda_\pm=\frac{v\pm\sqrt{v^2+4(\tu+1)}}2
\end{equation}
in a small neighborhood of the origin. Thus, 
 \eqref{KS2} is strictly hyperbolic. Then we can apply Lemma 2.1 in \cite{D15} to \eqref{KS2}, \eqref{ID2} to have the local existence result as follows.
 
 \begin{lemma} \label{lemma}
There are constants $\hat \delta>0$ and $\kappa>1$ such that whenever
\beq\label{sec2_1}
TV v_0(\cdot)+TV \tilde u_0(\cdot)< \hat\delta,
\eeq
the Cauchy problem \eqref{KS2}, \eqref{ID2} possesses a unique admissible BV solution $(v,\tilde u)$, defined on the strip $(-\infty,\infty)\times[0,1)$, taking values in a $\rho_0$-neighborhood of the origin and 
\beq\label{sec2_2}
TV v(\cdot,t)+TV \tilde u(\cdot,t)\leq \kappa[TV v_0(\cdot)+TV \tilde u_0(\cdot)],\qquad t\in[0,1).
\eeq
\end{lemma}

To extend the solution from local-in-time to global-in-time, we only need to study an admissible BV solution $(v,\tilde u)$ defined on $(-\infty,\infty)\times[0,T)$ for $T>0$. We want to show that the domain of $(v,\tilde u)$ can be extended to $(-\infty,\infty)\times[0,T+1)$. In the following sections we carry out analysis for the extension. As byproducts we also obtain the estimates \eqref{main_thm_1} and \eqref{main_thm_1a}. Theorem \ref{main_thm} is then proved.

From now on we assume that $(v,\tilde u)$ is an admissible BV solution of \eqref{KS2}, \eqref{ID2}, defined on a strip $(-\infty,\infty)\times[0,T)$ for some $T>0$, taking values in a $\rho_0$-neighborhood of the origin.

Our first step is to establish energy estimates assuming that the initial data satisfy \eqref{IDint}.
It starts with the entropy inequality for an admissible solution. It is straightforward to  verify that  \eqref{KS2} is endowed with a companion balance law
$$
\big[\frac{1}{2}v^2+(1+\tu)\ln(1+\tu)-\tu \big]_t+\big[v(1+\tu)\ln(1+\tu)\big]_x
+\tu(1+\tu)\ln(1+\tu)= 0.
$$
Here the convex entropy
$$
\eta=\frac{1}{2}v^2+(1+\tu)\ln(1+\tu)-\tu 
$$
was introduced in \cite{1d3} for the model with a diffusive but non-growth cellular population. Following the classical hyperbolic theory \cite{Smoller}, for small, admissible BV solutions, in which shock waves satisfy the Lax entropy condition, the following
entropy inequality holds on $(-\infty,\infty)\times[0,T)$:
\beq\label{en}
\big[\frac{1}{2}v^2+(1+\tu)\ln(1+\tu)-\tu \big]_t+\big[v(1+\tu)\ln(1+\tu)\big]_x
+\tu(1+\tu)\ln(1+\tu)\leq 0.
\eeq

Integrating \eqref{en} on $(-\infty,\infty)\times [0,t]$, $0\leq t< T$, we have 
\beq\label{en2}
\begin{split}
&\int_{-\infty}^\infty\big[\frac{1}{2}v^2+(1+\tu)\ln(1+\tu)-\tu \big](x,t)dx
+\int_0^t\int_{-\infty}^\infty \big[\tu(1+\tu)\ln(1+\tu)\big](x,\tau)dxd\tau\\
\leq&\int_{-\infty}^\infty\big[\frac{1}{2}v^2+(1+\tu)\ln(1+\tu)-\tu \big](x,0)dx.
\end{split}
\eeq

Note that by Taylor expansion,
\beq\label{en3}
(1+\tu)\ln(1+\tu)=\tu+\frac{1}{2}\tu^2+O(\tu^3), \qquad |\tu|< \rho_0.
\eeq
Thus, \eqref{en2} implies, for $0\le t<T$,
\beq\label{en4}
\begin{split}
\int_{-\infty}^\infty\big[v^2+\tu^2 \big](x,t)dx
&\leq C\int_{-\infty}^\infty\big[v_0^2+\tu_0^2 \big](x)dx\leq C\sigma^2,\\
\int_{0}^t\int_{-\infty}^\infty \tu^2(x,\tau) dxd\tau &\leq C\sigma^2,
\end{split}
\eeq
where $\sigma$ is defined in \eqref{IDint}.
Here and below, we use $C$ 
for a generic positive constant. In particular, $C$ is independent of $T$, $\rho_0$, $\sigma$ and $\delta$ (in \eqref{IDTV}).

Our next step is to construct a smooth, time-asymptotic solution to \eqref{KS2}, \eqref{ID2} by Chapman-Enskog expansion. Noting that $v$ is conserved and with nonzero mass, one does not expect time decay of its $L^1$ norm. It is necessary to construct a time-asymptotic solution with the $v$-component carrying the same mass, and extract it from $v$ to obtain the needed $L^1$-decay.

For this we identify leading terms in time decay rates in $\eqref{KS2}_2$ as
$$
v_x\approx-\tilde u.
$$
Therefore, with $\eqref{KS2}_1$ we define a time asymptotic  solution for $(v,\tilde u)$ as $(\bar v,\bar u)$, where $(\bar v,\bar u)$ satisfies
\begin{equation}\label{2.29}
\begin{cases}
\bar v_t+\bar u_x=0,\\
\bar v_x=-\bar u.
\end{cases}
\end{equation}
Substituting $\eqref{2.29}_2$ into $\eqref{2.29}_1$ gives us
\begin{equation}\label{2.30}
\begin{cases}
\bar v_t=\bar v_{xx},\\
\bar u=-\bar v_x.
\end{cases}
\end{equation}
Now we define $\bar v$ as the self-similar solution of $\eqref{2.30}_1$ carrying the same mass as $v_0$:
\beq\label{thdef}
\bar v(x,t)=\theta(x,t)=\frac{M}{\sqrt{4\pi(t+1)}} e^{-\frac{x^2}{4(t+1)}},\qquad
M=\int_{-\infty}^\infty v_0(x)\,dx,
\eeq
and by $\eqref{2.30}_2$,
$$
\bar u(x,t)=-\theta_x(x,t).
$$
Thus, we take $(\theta,-\theta_x)(x,t)$ as a time asymptotic solution to \eqref{KS2}, \eqref{ID2}, or equivalently, $(\theta,1-\theta_x)(x,t)$ as a time asymptotic solution to \eqref{KS}, \eqref{ID}.

The heat kernel $\theta$ defined in \eqref{thdef} plays a key role in our analysis. Note that 
\beq\label{Mes}
|M|\leq [\int_{-\infty}^\infty (1+x^2)^{-1} dx]^\frac{1}{2}
 [\int_{-\infty}^\infty (1+x^2)v_0^2(x) dx]^\frac{1}{2}\leq C\sigma
\eeq
by \eqref{IDint}.
Noting both $v$ and $\theta$ are conserved quantities carrying the same mass $M$, we introduce new variables so that the $v$-component is of zero mass.

Let
\beq\label{defw}
W=\left(
\begin{array}{c}
w_1\\
w_2
\end{array}
\right)
=
\left(
\begin{array}{c}
v-\theta\\
\tu
\end{array}
\right)
=
\left(
\begin{array}{c}
v-\theta\\
u-1
\end{array}
\right).
\eeq
By $\eqref{KS2}_1$ and \eqref{thdef}, we have
\beq\label{w1}
\int_{-\infty}^\infty w_1(x,t)\,dx=0.
\eeq

Under the new variables in \eqref{defw}, we rewrite \eqref{KS}-\eqref{ID}, or \eqref{KS2}-\eqref{ID2}, as
\beq\label{Weq}
\left\{
\begin{split}
&{w_1}_t+{w_2}_x +\theta_t =0,\\
&{w_2}_t+\big[(w_1+\theta)(1+w_2)]_x=-(1+w_2)w_2,
\end{split}
\right.
\eeq
with initial data
\beq\label{WID}
W_0(x)\equiv W(x,0)=(w_1,w_2)(x,0)=(w_{10},w_{20})(x)=(v_0(x)-\theta(x,0), u_0(x)-1).
\eeq

Still by \eqref{IDint}, we have
\[
\int_{-\infty}^\infty (1+x^2)\Big(
\big(w_{10}+\theta(x,0)\big)^2(x)+w_{20}^2(x)\Big)\, dx=\sigma^2.
\]
With \eqref{thdef} and \eqref{Mes}, it implies
\beq\label{2.15+}
\int_{-\infty}^\infty (1+x^2)
\Big(w_{10}^2(x)+w_{20}^2(x)\Big)\, dx\le C\sigma^2.
\eeq

\section{Weighted $L^2$ estimates and $L^1$ estimate for time decay}

In this section we prove several lemmas to estimate $w_1$ and $w_2$. 
First, we give the relative entropy inequality. By entropy inequality \eqref{en} and \eqref{Weq}, we can verify the following relative entropy inequality,
\beq\label{re0}
\begin{split}
&\big[\frac{1}{2}w_1^2+(1+w_2)\ln(1+w_2)-w_2+\theta_x w_2 \big]_t\\
+&\big[(\theta+w_1)(1+w_2)\ln(1+w_2)-\theta w_2+\theta_x w_1 +\theta_x(\theta+w_1)w_2\big]_x\\
+&w_2(1+w_2)\ln(1+w_2)+2\theta_x w_2 +\theta_x^2\\
\leq& \theta_{xt}w_2-\theta_x w_2^2+\theta_{xx}(\theta+w_1)w_2.
\end{split}
\eeq

\begin{lemma} \label{lemma3.1}
\beq\label{le1}
\int_{-\infty}^\infty (x^2+t+1)\big(
w^2_{1}(x,t)+w_{2}^2(x,t)\big)\, dx\leq C\sigma^2+CY,
\quad 0\leq t< T,
\eeq

\beq\label{le2}
\int_0^T \int_{-\infty}^\infty (x^2+\tau+1)(w_2+\theta_x)^2(x,\tau)\, dxd\tau\leq C\sigma^2+CY,
\eeq
where
\beq\label{Ydef}
Y=\int_0^T \int_{-\infty}^\infty w^2_{1}(x,\tau) dx d\tau.
\eeq
\end{lemma}
{\bf Proof:} Multiplying \eqref{re0} by $t$ gives us
\beq\label{re}
\begin{split}
&\left\{ t\big[\frac{1}{2}w_1^2+(1+w_2)\ln(1+w_2)-w_2+\theta_x w_2 \big]
\right\}_t
\\
+&\left\{ t\big[(\theta+w_1)(1+w_2)\ln(1+w_2)-\theta w_2+\theta_x w_1 +\theta_x(\theta+w_1)w_2\big]\right\}_x\\
+&t\left\{w_2(1+w_2)\ln(1+w_2)+2\theta_x w_2 +\theta_x^2\right\}\\
\leq& R,
\end{split}
\eeq
where 
\[
R=\frac{1}{2}w_1^2+(1+w_2)\ln(1+w_2)-w_2+\theta_x w_2
+t\big((\theta_{xt}+\theta_{xx}\theta)w_2-\theta_x w_2^2+\theta_{xx}w_1w_2\big).
\]
Integrating \eqref{re} on $(-\infty,\infty)\times [0,t],$ $0\leq t<T$, we have 
\beq\label{re2}
\begin{split}
&\int_{-\infty}^\infty t\big[\frac{1}{2}w_1^2+(1+w_2)\ln(1+w_2)-w_2 \big](x,t)dx
\\
+&\int_{0}^t\int_{-\infty}^\infty\tau \big[w_2(1+w_2)\ln(1+w_2)+2\theta_x w_2 +\theta_x^2\big](x,\tau)dxd\tau\\
\leq& \int_{0}^t\int_{-\infty}^\infty R(x,\tau)\, dxd\tau-\int_{-\infty}^\infty t (\theta_x w_2)(x,t)\, dx.
\end{split}
\eeq

We estimate each term in \eqref{re2} as follows. Applying \eqref{en3}, we have
\beq\label{25}
\int_{-\infty}^\infty t\big[\frac{1}{2}w_1^2+(1+w_2)\ln(1+w_2)-w_2 \big](x,t)dx\geq 
\int_{-\infty}^\infty t(\frac{1}{2}w_1^2+\frac{1}{4}w_2^2 )(x,t)dx,
\eeq
where we have used the assumption $|w_2|< \rho_0$ for $\rho_0$ sufficiently small. We also have
\beq\label{26}
\begin{split}
&\int_{0}^t\int_{-\infty}^\infty\tau \big[w_2(1+w_2)\ln(1+w_2)+2\theta_x w_2 +\theta_x^2\big](x,\tau)dxd\tau\\
\ge\,&\int_{0}^t\int_{-\infty}^\infty\tau (w_2+\theta_x)^2(x,\tau)dxd\tau-C\int_{0}^t\int_{-\infty}^\infty\tau |w_2(x,\tau)|^3dxd\tau\\
\ge\,&\int_{0}^t\int_{-\infty}^\infty\tau (w_2+\theta_x)^2(x,\tau)dxd\tau-C\int_{0}^t\int_{-\infty}^\infty\tau |w_2+\theta_x|^3(x,\tau)dxd\tau\\
&-C\int_{0}^t\int_{-\infty}^\infty\tau |\theta_x|^3(x,\tau)dxd\tau\\
\ge\,&[1-C(\rho+|M|)]\int_{0}^t\int_{-\infty}^\infty\tau (w_2+\theta_x)^2(x,\tau)dxd\tau-C|M|^3\\
\ge\,&
\frac{1}{2}\int_{0}^t\int_{-\infty}^\infty \tau(w_2+\theta_x)^2(x,\tau)dxd\tau-C\sigma^3,
\end{split}
\eeq
where we have applied Young's inequality, \eqref{thdef} and \eqref{Mes}, and $\sigma$ in \eqref{IDint} is sufficiently small.

Then we estimate the right hand side of \eqref{re2}. Using \eqref{thdef}, \eqref{Mes} and \eqref{en4}, we have 
\beq
\begin{split}\label{27}
&|-\int_{-\infty}^\infty t(\theta_x w_2)(x,t)dx|\\
\leq& C\sigma \int_{-\infty}^\infty e^{-\frac{x^2}{8(t+1)}}|w_2(x,t)|dx\\
\leq& C\sigma^2 \int_{-\infty}^\infty (t+1)^{-1}e^{-\frac{x^2}{4(t+1)}}dx
+\frac{1}{8} \int_{-\infty}^\infty (t+1) w^2_2(x,t)dx\\
\leq &C \sigma^2+\frac{1}{8} \int_{-\infty}^\infty tw^2_2(x,t)dx.
\end{split}
\eeq

We divide $R$ into parts and consider the following items.
First, using \eqref{Ydef}, \eqref{en4}, \eqref{thdef} and \eqref{Mes}, we have
\beq
\begin{split}\label{28}
&\int_{0}^t\int_{-\infty}^\infty |\frac{1}{2}w_1^2+(1+w_2)\ln(1+w_2)-w_2+\theta_x w_2|(x,\tau)\,dxd\tau\\
\leq&\int_{0}^t\int_{-\infty}^\infty (\frac{1}{2}w_1^2+\frac{3}{2}w_2^2
+\frac{1}{2}\theta_x^2)(x,\tau)\,dxd\tau\\
\leq& \frac{1}{2}Y+C\sigma^2,
\end{split}
\eeq
when $0\leq t\leq T$.

Secondly, still using \eqref{en4}, \eqref{thdef} and \eqref{Mes} we have
\beq\label{29}
\begin{split}
&\int_{0}^t\int_{-\infty}^\infty
\tau|(\theta_{x\tau}+\theta_{xx}\theta)w_2|(x,\tau)\,dxd\tau\\
\leq&\frac{1}{2}\int_{0}^t\int_{-\infty}^\infty
\tau^2 (\theta_{x\tau}+\theta_{xx}\theta)^2(x,\tau)\,dxd\tau
+
\frac{1}{2}\int_{0}^t\int_{-\infty}^\infty w^2_2(x,\tau)\,dxd\tau\\
\leq& C\sigma^2,
\end{split}
\eeq

\beq\label{30}
\int_{0}^t\int_{-\infty}^\infty
\tau|\theta_{x}w_2^2|(x,\tau)\,dxd\tau
\leq C\sigma\int_{0}^t\int_{-\infty}^\infty w_2^2(x,\tau)\,dxd\tau
\leq C\sigma^3,
\eeq
and
\beq\label{31}
\begin{split}
&\int_{0}^t\int_{-\infty}^\infty
\tau|\theta_{xx}w_1 w_2|(x,\tau)\,dxd\tau\\
\leq&C\sigma\int_{0}^t\int_{-\infty}^\infty (\tau+1)^{-\frac{1}{2}}|w_1 w_2|(x,\tau)\,dxd\tau\\
\leq&C\sigma\int_{0}^t\int_{-\infty}^\infty \frac{1}{2}(w^2_1+ w^2_2)(x,\tau)\,dxd\tau\\
\leq& C\sigma Y+C\sigma^3.
\end{split}
\eeq

Combining these estimates together, we have
\beq\label{est}
\int_{-\infty}^\infty t(w_1^2+w_2^2 )(x,t)dx+\int_{0}^t\int_{-\infty}^\infty \tau(w_2+\theta_x)^2(x,\tau)dxd\tau
\leq CY+C\sigma^2.
\eeq

Similarly, we multiply \eqref{re0} by $x^2$ to have
\beq\label{33}
\begin{split}
&\Big\{ x^2\big[\frac{1}{2}w_1^2+(1+w_2)\ln(1+w_2)-w_2+\theta_x w_2 \big]\Big\}_t\\
+&\Big\{ x^2\big[(\theta+w_1)(1+w_2)\ln(1+w_2)-\theta w_2+\theta_x w_1 +\theta_x(\theta+w_1)w_2\big]\Big\}_x\\
+&x^2\left\{w_2(1+w_2)\ln(1+w_2)+2\theta_x w_2 +\theta_x^2\right\}\\
\leq& S,
\end{split}
\eeq
where 
\[
S=2x\big[(\theta+w_1)(1+w_2)\ln(1+w_2)-\theta w_2+\theta_x w_1 +\theta_x(\theta+w_1)w_2\big]
+x^2\big((\theta_{xt}+\theta_{xx}\theta)w_2-\theta_x w_2^2+\theta_{xx}w_1w_2\big).
\]
Integrating \eqref{33} on $(-\infty,\infty)\times [0,t],$ $0\leq t< T$, we have 
\beq\label{34}
\begin{split}
&\int_{-\infty}^\infty x^2\big[\frac{1}{2}w_1^2+(1+w_2)\ln(1+w_2)-w_2+\theta_x w_2 \big](x,t)\, dx
\\
+&\int_{0}^t\int_{-\infty}^\infty x^2\left\{w_2(1+w_2)\ln(1+w_2)+2\theta_x w_2 +\theta_x^2\right\}(x,\tau)\, dx d\tau\\
\leq& \int_{0}^t\int_{-\infty}^\infty S(x,\tau)\, dxd\tau+\int_{-\infty}^\infty x^2 \Big[\frac{1}{2}w_1^2+(1+w_2)\ln(1+w_2)-w_2+\theta_x w_2\Big](x,0)\, dx.
\end{split}
\eeq

Similar to \eqref{25} and \eqref{26} we have
\beq\label{35}
\int_{-\infty}^\infty x^2\big[\frac{1}{2}w_1^2+(1+w_2)\ln(1+w_2)-w_2 \big](x,t)dx\geq 
\int_{-\infty}^\infty x^2(\frac{1}{2}w_1^2+\frac{1}{4}w_2^2 )(x,t)dx,
\eeq
\beq\label{36}
\begin{split}
&\int_{0}^t\int_{-\infty}^\infty x^2 \big[w_2(1+w_2)\ln(1+w_2)+2\theta_x w_2 +\theta_x^2\big](x,\tau)dxd\tau\\
\ge\,&
\frac{1}{2}\int_{0}^t\int_{-\infty}^\infty x^2(w_2+\theta_x)^2(x,\tau)dxd\tau-C\sigma^3,
\end{split}
\eeq
\beq\label{37}
\begin{split}
&\int_{-\infty}^\infty x^2\big[\frac{1}{2}w_1^2+(1+w_2)\ln(1+w_2)-w_2+\theta_xw_2 \big](x,0)dx\\
\le \,&
C\int_{-\infty}^\infty x^2\big[w_{10}^2(x)+w_{20}^2(x)+\theta_x^2(x,0)\big]dx\\
\le \,&C\sigma^2,
\end{split}
\eeq
where we have applied \eqref{2.15+} in the last estimate.
Similar to \eqref{27} we have
\beq
\begin{split}\label{38}
&|-\int_{-\infty}^\infty x^2(\theta_x w_2)(x,t)dx|\\
\leq& C\sigma^2 \int_{-\infty}^\infty \frac{x^2}{(t+1)^2}e^{-\frac{x^2}{4(t+1)}}dx+\frac18\int_{-\infty}^\infty x^2w_2^2(x,t)dx\\
\leq& C\sigma^2+\frac18\int_{-\infty}^\infty x^2w_2^2(x,t)dx.
\end{split}
\eeq

For the double integral on the right-hand side of \eqref{34} we consider each part of $S$ as follows. Applying \eqref{en3}, \eqref{en4}, \eqref{thdef}, \eqref{Mes} and Young's inequality, we have
\beq
\begin{split}\label{39}
&\int_0^t\int_{-\infty}^\infty\Big|2x\big[(\theta+w_1)(1+w_2)\ln(1+w_2)-\theta w_2+\theta_x w_1 +\theta_x(\theta+w_1)w_2\big]\Big|(x,\tau)dxd\tau\\
=\,&\int_0^t\int_{-\infty}^\infty2|x|\big|(\theta+w_1)[w_2+O(w_2^2)]-\theta w_2+\theta_x w_1 +\theta_x(\theta+w_1)w_2\big|(x,\tau)dxd\tau\\
=\,&\int_0^t\int_{-\infty}^\infty2|x|\big[|w_1(w_2+\theta_x)|+|(\theta+w_1)O(w_2^2)|+|\theta\theta_x w_2|+|\theta_xw_1w_2|\big](x,\tau)dxd\tau\\
\le\,&8\int_0^t\int_{-\infty}^\infty w_1^2(x,\tau)dxd\tau+\frac18\int_0^t\int_{-\infty}^\infty x^2(w_2+\theta_x)^2(x,\tau)dxd\tau\\
&+C\int_0^t\int_{-\infty}^\infty|x|( |w_1|w_2^2)(x,\tau)dxd\tau+C\int_0^t\int_{-\infty}^\infty(w_1^2+w_2^2+\theta_x^2)(x,\tau)dxd\tau\\
\le\,&8Y+\frac18\int_0^t\int_{-\infty}^\infty x^2(w_2+\theta_x)^2(x,\tau)dxd\tau+C\int_0^t\int_{-\infty}^\infty|x|\big|w_1w_2(w_2+\theta_x)\big|(x,\tau)dxd\tau\\
&+C\int_0^t\int_{-\infty}^\infty(w_1^2+w_2^2+\theta_x^2)(x,\tau)dxd\tau\\
\le\,&8Y+\frac14\int_0^t\int_{-\infty}^\infty x^2(w_2+\theta_x)^2(x,\tau)dxd\tau+C\int_0^t\int_{-\infty}^\infty(w_1^2+w_2^2+\theta_x^2)(x,\tau)dxd\tau\\
\le\,&CY+\frac14\int_0^t\int_{-\infty}^\infty x^2(w_2+\theta_x)^2(x,\tau)dxd\tau+C\sigma^2.
\end{split}
\eeq
Similar to \eqref{29}-\eqref{31} we also have
\beq
\begin{split}\label{40}
&\int_0^t\int_{-\infty}^\infty x^2\big|(\theta_{xt}+\theta_{xx}\theta)w_2-\theta_x w_2^2+\theta_{xx}w_1w_2\big|
(x,\tau)dxd\tau\\
\le\,&\frac12\int_0^t\int_{-\infty}^\infty x^4(\theta_{xt}+\theta_{xx}\theta)^2(x,\tau)dxd\tau+C\int_0^t\int_{-\infty}^\infty(w_1^2+w_2^2)(x,\tau)dxd\tau\\
\le\,&C\sigma^2+CY.
\end{split}
\eeq

Substituting \eqref{35}-\eqref{40} into \eqref{34} and simplify, we arrive at
\beq\label{est2}
\int_{-\infty}^\infty x^2(w_1^2+w_2^2 )(x,t)dx+\int_{0}^t\int_{-\infty}^\infty x^2(w_2+\theta_x)^2(x,\tau)dxd\tau
\leq CY+C\sigma^2.
\eeq

From \eqref{en4} and \eqref{thdef}-\eqref{defw}, we also have 
\beq\label{est3}
\int_{-\infty}^\infty (w_1^2+w_2^2 )(x,t)dx+\int_{0}^t\int_{-\infty}^\infty (w_2+\theta_x)^2(x,\tau)dxd\tau
\leq C\sigma^2.
\eeq

Finally, the lemma is proved by combining \eqref{est2}, \eqref{est} and \eqref{est3}. $\square$

\begin{lemma}\label{lemma3.2}
\beq\label{Yest}
Y\leq C\sigma^2.
\eeq
\end{lemma}

\pf
We introduce a potential function 
\beq\label{Phi}
\Psi(x,t)=\int_{-\infty}^x w_1(y,t)\, dy.
\eeq
So by \eqref{w1}, we know that $\Psi(\pm\infty,t)=0$. From \eqref{Weq} we have
\[
\Psi_t=-w_2-\theta_x,\qquad
\Psi_x=w_1.
\]
Therefore,
\[
\begin{split}
&(\frac{1}{2}\Psi^2-w_2\Psi)_t
\\
=&\Psi(\Psi_t-{w_2}_t)-w_2\Psi_t\\
=&\Psi\Big(-w_2-\theta_x+\big((w_1+\theta)(1+w_2)\big)_x+(1+w_2)w_2\Big)+w^2_2+w_2\theta_x\\
=&\Psi\Big(\big(w_1+w_1w_2+\theta w_2\big)_x+w^2_2\Big)+w^2_2+w_2\theta_x\\
=&\big(\Psi(w_1+w_1w_2+\theta w_2)\big)_x-w_1(w_1+w_1w_2+\theta w_2)+(\Psi+1)w_2^2+w_2\theta_x,
\end{split}
\]
which gives us
$$
w_1(w_1+w_1w_2+\theta w_2)=(w_2\Psi-\frac{1}{2}\Psi^2)_t+\big(\Psi(w_1+w_1w_2+\theta w_2)\big)_x
+w_2\theta_x+(\Psi+1)w_2^2.
$$
Integrating the equation on $(-\infty,\infty)\times [0,t],$ $0\leq t<T$, we arrive at 
\beq
\begin{split}\label{3.64}
&\int_0^t\int_{-\infty}^\infty [w_1(w_1+w_1w_2+\theta w_2)](x,\tau)dxd\tau\\
=\,&\int_{-\infty}^\infty(w_2\Psi-\frac{1}{2}\Psi^2)(x,t)dx+\int_0^t\int_{-\infty}^\infty(w_2\theta_x)(x,\tau)dxd\tau\\
&+\int_0^t\int_{-\infty}^\infty[(\Psi+1)w_2^2](x,\tau)dxd\tau-\int_{-\infty}^\infty(w_2\Psi-\frac{1}{2}\Psi^2)(x,0)dx.
\end{split}
\eeq

For the left hand side of \eqref{3.64}, by choosing $\rho_0<\frac{1}{2}$ and hence $|w_2|< \frac{1}{2}$, we have
\beq
\begin{split}\label{3.65}
&\int_0^t\int_{-\infty}^\infty [w_1(w_1+w_1w_2+\theta w_2)](x,\tau)dxd\tau \\ \geq& \int_0^t\int_{-\infty}^\infty (\frac{1}{2}w_1^2+\theta w_1 w_2)(x,\tau)dxd\tau \\
\geq& \int_0^t\int_{-\infty}^\infty (\frac{1}{4}w_1^2-\theta^2 w_2^2)(x,\tau)dxd\tau \\
\geq &
\frac{1}{4} \int_0^t\int_{-\infty}^\infty w_1^2(x,\tau)dxd\tau-C\sigma^2,
\end{split}
\eeq
where we have applied \eqref{en4} and \eqref{thdef}.

On the right-hand side of \eqref{3.64}, also by \eqref{en4} we have
\beq
\label{3.66}
\int_{-\infty}^\infty(w_2\Psi-\frac{1}{2}\Psi^2)(x,t)dx\le \frac{1}{2}\int_{-\infty}^\infty w_2^2(x,t) dx\le C\sigma^2.
\eeq
It is straightforward to verify that 
\beq
\label{3.67}
|\int_0^t\int_{-\infty}^\infty
(w_2\theta_x)(x,\tau)
dxd\tau|\leq
(\int_0^t\int_{-\infty}^\infty
w_2^2(x,\tau)
dxd\tau)^{\frac{1}{2}}
(\int_0^t\int_{-\infty}^\infty
\theta_x^2(x,\tau)
dxd\tau)^{\frac{1}{2}}
\leq
C\sigma^2.
\eeq

Next, we give a bound on $\Psi$. From \eqref{le1},
\beq\label{arg}
(\int_{-\infty}^\infty |w_1(x,t)|dx)^2
\leq \int_{-\infty}^\infty (x^2+t+1)^{-1}dx \int_{-\infty}^\infty(x^2+t+1) w^2_1(x,t)\, dx\leq C(\sigma^2+Y).
\eeq
Therefore,
\[
|\Psi(x,t)|\le \int_{-\infty}^\infty |w_1(x,t)|dx\leq C(\sigma+\sqrt Y).
\]
By \eqref{en4}, we further have
\beq\label{3.69}
|\int_0^t \int_{-\infty}^\infty[(\Psi+1)w_2^2](x,\tau) dxd\tau  |\leq C(\sigma+\sqrt  Y+1)\int_0^t \int_{-\infty}^\infty w_2^2(x,\tau) dxd\tau \leq C\sigma^2(1+\sqrt Y).
\eeq

Finally, we need to prove that 
\beq\label{phi2}
\int_{-\infty}^\infty\Psi^2(x,0)dx \leq C \sigma^2.
\eeq

First, by integration by parts and noting that $v_0$ and hence $w_{10}$ decays to zero sufficiently fast as $x\rightarrow-\infty$ as assumed in Theorem \ref{main_thm}, we have
\begin{equation}\label{3.40+}
\begin{split}
&\int^0_{-\infty}[\int_{-\infty}^x w_1(y,0)dy]^2 dx
=-2\int_{-\infty}^0xw_1(x,0)\int_{-\infty}^xw_1(y,0)\,dydx\\
\leq &2\int^0_{-\infty} x^2w_1^2(x,0)\,dx+\frac12\int^0_{-\infty}[\int_{-\infty}^x w_1(y,0)dy]^2 dx.
\end{split}
\end{equation}
Equation \eqref{3.40+} is simplified to
\begin{equation}\label{3.40++}
\int^0_{-\infty}[\int_{-\infty}^x w_1(y,0)dy]^2 dx
\leq 4 \int^0_{-\infty} x^2
w^2_1(x,0)dx\leq C\sigma^2
\end{equation}
by \eqref{2.15+}.  

 Recall \eqref{w1},
\[
\int_{-\infty}^\infty w_1(y,t)\,dy=0.
\]
Then, we rewrite $\Psi$ in \eqref{Phi} as
$$
\Psi(x,t)=-\int_x^\infty w_1(y,t)\,dy.
$$
A bound similar to \eqref{3.40++} holds for the integral of $\Psi^2(x,0)$ in $(0,\infty)$, and we prove \eqref{phi2}.

In summary, substituting \eqref{3.65}-\eqref{3.67} and \eqref{3.69}-\eqref{phi2} into \eqref{3.64} and simplifying, we arrive at
\[
Y\leq C\sigma^2+C\sigma^2 \sqrt Y\le C\sigma^2+C\sigma^4 + \frac12Y,
\]
which is simplified to 
\[
Y\leq C\sigma^2
\]
for some $C$.  $\square$

Combining Lemmas \ref{lemma3.1} and \ref{lemma3.2}, we upgrade \eqref{le1} and \eqref{le2} to
\beq\label{le1n}
\int_{-\infty}^\infty (x^2+t+1)\big(
w^2_{1}(x,t)+w_{2}^2(x,t)\big)\, dx\leq C\sigma^2,
\quad 0\leq t< T,
\eeq

\beq\label{le2n}
\int_0^T \int_{-\infty}^\infty (x^2+\tau+1)(w_2+\theta_x)^2(x,\tau)\, dxd\tau\leq C\sigma^2.
\eeq

Finally, we can prove an $L^1$ estimate on $(w_1,w_2)$.
\begin{lemma}\label{lemma3.3} 
\beq\label{l1w}
\int_{-\infty}^\infty
\Big(|w_{1}(x,t)|+|w_{2}(x,t)|\Big)\, dx\leq C\sigma(t+1)^{-\frac{1}{4}}, \qquad
\quad 0\leq t< T,
\eeq
i.e.
\beq\label{l1uv}
\int_{-\infty}^\infty
\Big(|v(x,t)-\theta(x,t)|+|\tilde u(x,t)|\Big)\, dx\leq C\sigma(t+1)^{-\frac{1}{4}} ,\qquad
\quad 0\leq t< T,
\eeq
using the relation \eqref{defw}, where $\theta$ is defined in \eqref{thdef}. 
\end{lemma}
Proof: In fact, applying Cauchy–Schwarz inequality and \eqref{le1n} we have
\[
\begin{split}
&\bigg[\int_{-\infty}^\infty \big(|w_1(x,t)|+|w_2(x,t)|\big)dx\bigg]^2\\
\leq& 2\int_{-\infty}^\infty (x^2+t+1)^{-1}dx \int_{-\infty}^\infty(x^2+t+1) (w^2_1(x,t)+w^2_2(x,t))dx\\
\leq& C(t+1)^{-\frac{1}{2}}\sigma^2.
\end{split}
\]
Therefore,
\[
\int_{-\infty}^\infty \big(|w_1(x,t)|+|w_2(x,t)|\big)dx\leq C
\sigma(t+1)^{-\frac{1}{4}}.
\]

\section{Bounds on the Total Variation}

We continue to work on the admissible BV solution $(v, \tu)$ of \eqref{KS2}, \eqref{ID2}, defined on $(-\infty, \infty) \times [0, T )$ for $T > 0$ and taking values in a $\rho_0$-neighborhood of the origin.  The solution satisfies estimates obtained in Section 3, in particular \eqref{l1w}, with $W$ defined in \eqref{defw}. In addition, we assume for the moment that
\begin{equation}\label{4.1}
TV v(\cdot,t)+TV \tilde u(\cdot,t)\leq \varepsilon_0,\qquad t\in[0,T),
\end{equation}
or equivalently,
\begin{equation}\label{4.2}
TV w_1(\cdot,t)+TV w_2(\cdot,t)\leq \varepsilon_1,\qquad t\in[0,T),
\end{equation}
where $\varepsilon_0$ and $\varepsilon_1$ are some sufficiently small numbers. The value of $\varepsilon_1$ will be determined when we go through this section, and \eqref{4.2} will be justified at the end of the section.

Our goal is to obtain needed estimates on the total variation for the extension from local solutions to global solutions, and hence prove Theorem \ref{main_thm}. While we generally follow the road map from \cite{D06}, our focus is on new details related to the nonzero-mass situation. We also incorporate the $L^1$ decay estimate \eqref{l1w} to upgrade bounds on the total variation to achieve \eqref{main_thm_1}, which is not available in \cite{D06,D13,D15}. For details that are completely parallel to those in \cite{D06}, we give a brief outline, and readers are referred to it.

\subsection{Redistribution of Dissipation}
Noting that the first equation of \eqref{KS2} (or equivalently, of \eqref{Weq}) has no explicit dissipation, we apply Dafermos' idea \cite{D06,D13,D14,D15} to re-distribute the uneven dissipation between its two equations.
This is done by introducing a change of variables through a nonlocal transformation via the potential function $\Psi$ defined in \eqref{Phi}. Let
\begin{equation}\label{4.3}
\phi(x,t)=\frac{1}{2}\Psi(x,t)=\frac{1}{2}\int^x_{-\infty} w_1(y,t) dy, \quad -\infty<x<\infty,\quad 0\leq t<T.
\end{equation}
Correspondingly, we have  derivatives
\begin{equation}\label{4.4}
\phi_t=-\frac{1}{2}w_2-\frac{1}{2}\theta_x,\qquad
\phi_x=\frac{1}{2}w_1,
\end{equation}
where we have applied \eqref{Weq} and \eqref{heat}.
Now we set
\begin{equation}\label{4.5}
\Phi(x,t)=\left(
\begin{array}{c}
0\\
\phi(x,t)
\end{array}
\right).
\end{equation}
Our new variable to replace $W$ is 
\beq
\label{Weq1}
\hat W(x,t)=W(x,t)-\Phi(x,t)=
\left(
\begin{array}{c}
w_1\\
w_2-\phi
\end{array}
\right)(x,t)
\equiv
\left(
\begin{array}{c}
w_1\\
w_3
\end{array}
\right)(x,t).
\eeq

Using \eqref{defw} we rewrite \eqref{Weq} in vector notations as
\beq\label{4.7}
\partial_t W(x,t)+\partial_x F(W(x,t),x,t)+G(W(x,t),x,t)=0,
\eeq
with 
\begin{equation}\label{4.8}
F(W,x,t)=
\left(
\begin{array}{c}
w_2\\
(w_1+\theta)(1+w_2)
\end{array}
\right),
\quad
G(W,x,t)=
\left(
\begin{array}{c}
\theta_t\\
(1+w_2)w_2
\end{array}
\right),
\end{equation}
where $F$ and $G$ depend on $x$ and $t$ explicitly through $\theta(x,t)$ defined in \eqref{main_thm_1b}.

Under the new variable $\hat W$ we write \eqref{4.7} as
\beq\label{Weq2}
\partial_t \hat W(x,t)+\partial_x  \hat F( \hat W(x,t),\Phi(x,t),x,t)+  \hat G( \hat W(x,t),\Phi(x,t),x,t)=0,
\eeq
with 
\beq\label{Weq3}
 \hat F(\hat W,\Phi,x,t)=F(\hat W+\Phi,x,t)-F(\Phi,x,t),
\eeq
and
\beq\label{Weq4}
 \hat G(\hat W,\Phi,x,t)=
\left(
\begin{array}{c}
\frac{1}{2}w_1+\theta_{xx}\\
\frac12w_3+\frac12\theta w_1+\frac12\phi+\theta_x(\frac12+\phi)+(w_3+\phi)^2
\end{array}
\right).
\eeq
Here we have applied \eqref{heat} and \eqref{4.4}-\eqref{4.8}. The corresponding initial data are from \eqref{WID} as
\begin{equation}\label{4.12}
\hat W(x,0)=W_0(x)-\Phi(x,0).
\end{equation}

In contrast of \cite{D14}, $\hat F$ in \eqref{Weq3} contains $\Phi$, which does not have sufficient regularity to allow the application of results from \cite{DH82}. On the other hand, the detailed proof on $TV$-bounds in \cite{D06} is for the case where $\hat F$ does not depend on $x$ and $t$ explicitly. We now follow the general approach in \cite{D06} but provide new details specific to our $\hat F$ and upgrade $TV$ bounds with time-decay rates.

Our strategy is the following. We regard $\Phi$ as a known function of $x$ and $t$ via \eqref{4.3}, \eqref{4.5} and the known solution $W(x,t)$ of \eqref{4.7} , \eqref{WID} on $(-\infty,\infty)\times[0,T)$. With the known inhomogeneity through $\Phi$ and $\theta$ we reconstruct the solution $\hat W$ of \eqref{Weq2}, \eqref{4.12}. It is done by constructing a sequence of approximate solutions and then taking its limit. The total variations of the approximate solutions can be controlled and estimated, and hence we obtain bounds on the total variation of the limit. After showing that the limit is the solution of  \eqref{Weq2}, \eqref{4.12}, we obtain the needed bounds on the total variation of $W$ by the uniqueness and \eqref{Weq1}. Then we can extend the time interval in the local existence result, Lemma \ref{lemma}, and prove the time-decay rate in \eqref{main_thm_1}.

Next, we have some preparation for our analysis.
From \eqref{4.8} we find that
\begin{equation}\label{4.13}
F_{W}(W,x,t)=\begin{pmatrix}0&1\\1+w_2&w_1+\theta\end{pmatrix},\qquad
F_{W}^{-1}(W,x,t)=\frac1{1+w_2}\begin{pmatrix}-(w_1+\theta)&1\\1+w_2&0\end{pmatrix}.
\end{equation}
Here $F_{W}(W,x,t)$ has two real, distinct eigenvalues $\lambda_\pm(W,x,t)=\lambda_\pm$, where $\lambda_\pm$ are given in \eqref{2.1}. That is,
\begin{equation}\label{4.14}
\lambda_\pm(W,x,t)=\frac12[w_1+\theta\pm\sqrt{(w_1+\theta)^2+4(1+w_2)}].
\end{equation}
The corresponding eigenvectors are $r_\pm(W,x,t)=r_\pm=(1,\lambda_\pm)^t$, which form a matrix with $r_\pm$ as columns,
\begin{equation}\label{4.15}
R(W,x,t)=\begin{pmatrix}1&1\\ \lambda_-(W,x,t)&\lambda_+(W,x,t)\end{pmatrix}.
\end{equation}
It is straightforward to verify that $\nabla \lambda_\pm\cdot r_\pm\neq0$, and hence the two characteristic families of \eqref{4.7} are genuinely nonlinear. Similarly, \eqref{Weq2} is strictly hyperbolic, and its two  characteristic families are genuinely nonlinear. 

From \eqref{Weq4} we find
\begin{equation}\label{4.16}
\hat G_{\hat W}(\hat W,\Phi,x,t)=\begin{pmatrix}\frac12&0\\ \frac12\theta&\frac12+2(w_3+\phi)\end{pmatrix}.
\end{equation}
This is to be compared with $G_W(W,x,t)$, where the first row is zero, see \eqref{4.8}.
In contrast, \eqref{4.16} is diagonally dominant in the $\rho_0$-neighborhood of the origin. Thus, we have indeed redistributed the dissipation evenly between the two equations in \eqref{Weq2}.

From \eqref{4.16} and by similar calculation, we have
\begin{equation}\label{4.17}
\begin{split}
&\hat G_{\hat W}(\hat W,\Phi,x,t)=O(1),\\
&\hat G_{\Phi}(\hat W,\Phi,x,t)=O(1)(1+|\theta_x|+|w_3|+|\phi|)=O(1),\\
&\hat G_x(\hat W,\Phi,x,t)=O(1)(|\theta_{xxx}|+|\theta_x||\hat W|+|\theta_{xx}|+|\theta_{xx}||\phi|)=O(1).
\end{split}
\end{equation}

\subsection{Approximate Solutions}
Our approach to construct a sequence of approximate solutions is a combination of Glimm's random choice scheme \cite{Glimm} and the operator splitting techniques by Dafermos and Hsiao \cite{DH82} and Dafermos \cite{D06}.

We continue to use $C$ for a generic positive constant, independent of $T$, $\varepsilon_1$ (or $\varepsilon_0$) and initial data. For vectors we use 1-norm and for matrices we use the corresponding  induced norm. An approximate solution $\hat W_h$ is constructed as follows.

First, one generates a random sequence $\{\zeta_m\}$ of points in the interval $(-1,1)$. Fixing a spatial mesh-length $h$, which will serve as the (eventually vanishing) parameter, we select a temporal mesh-length $\tau=h/\lambda$, where $\lambda\ge2$ is a fixed positive number so that  waves emanating from points at a distance of $2h$ cannot collide on a time interval of length $\tau$.

We partition the upper half-plane into strips
\[
\mathcal A_m=\{(x,t): -\infty<x<\infty, \,m\tau< t< (m+1)\tau
\},\quad m=0,1, 2,\dots, 
\]
and identify the mesh-points $(kh,m\tau)$, for $k+m$ even, and the random mesh-points $(y_{k,m}, m\tau)$, for $k+m$ odd, where
\[
y_{k,m}=(k+\zeta_m)h.
\]

Let $m^*$ be the largest integer with $m^*\tau< T$. For $k=0,\pm1,\pm2,\dots$, and $m=0,1, \dots, m^*$, we set
\begin{equation}\label{4.18}
\Phi_{k,m}=\Phi(kh,m\tau).
\end{equation}

We begin the algorithm by setting 
\begin{equation}\label{4.19}
\hat W_h(x,0-)=\hat W_0(x)=W_0(x)-\Phi(x,0)=\begin{pmatrix}w_{10}(x)\\ w_{20}(x)-\phi(x,0)\end{pmatrix}, \qquad-\infty< x<\infty,
\end{equation}
see \eqref{WID}. Assuming now $\hat W_h$ has been determined on $\cup_{l=0}^{m-1}\mathcal A_l$, for some $m\le m^*$, we extend its domain to the next strip $\mathcal A_m$ as follows.

First, for any $k$ with $k+m$ odd, we set 
\begin{eqnarray}
&&\hat{W}_{k,m}=\hat{W}_h(y_{k,m},m\tau-),\qquad W_{k,m}=\hat{W}_{k,m}+\Phi_{k,m},\label{4.20}\\
&&S_{k,m}=F_W^{-1}(W_{k,m}, kh,m\tau)F_W(\Phi_{k,m},kh,m\tau),\label{4.21}\\
&&\tilde S_{k,m}=F_W^{-1}(W_{k,m}, kh,m\tau)[F_x(W_{k,m},kh,m\tau)-F_x(\Phi_{k,m},kh,m\tau)],\label{4.22}\\
&&G_{k,m}=\hat G(\hat W_{k,m}, \Phi_{k,m},kh,m\tau).\label{4.23}
\end{eqnarray}
Then we set
\begin{equation}\label{4.24}
\begin{cases}
W^L_{k,m}=W_{k,m}+S_{k,m}[\Phi_{k-1,m}-\Phi_{k,m}]
+h\tilde S_{k,m}-\tau G_{k,m},
\\
W^R_{k,m}=W_{k,m}+S_{k,m}[\Phi_{k+1,m}-\Phi_{k,m}]
-h\tilde S_{k,m}-\tau G_{k,m}.
\end{cases}
\end{equation}
Note that when comparing with \cite{D06}, the terms $\pm h\tilde S_{k,m}$ on the right-hand side are additional  to account for the explicit dependence of $F$ on $x$ and $t$ through $\theta$, a consequence of the nonzero-mass situation under consideration. Since our new details are closely related to those terms, we substitute \eqref{4.8} and \eqref{4.13} into \eqref{4.22} to have the explicit formulation
\begin{equation}\label{4.25}
\tilde S_{k,m}=\frac1{1+(w_2)_{k,m}}\begin{pmatrix}\theta_x(kh,m\tau)(\hat w_2)_{k,m}\\0\end{pmatrix}.
\end{equation}
Here $(w_2)_{k,m}$ stands for the second component of $W_{k,m}$, and $(\hat w_2)_{k,m}$ for the second component of $\hat W_{k,m}$, etc.

Now we denote 
\begin{equation}\label{4.26}
\begin{cases}
{\hat W}^L_{k,m}=W^L_{k,m}-\Phi_{k-1,m},
\\
{\hat W}^R_{k,m}=W^R_{k,m}-\Phi_{k+1,m}.
\end{cases}
\end{equation}
One can verify that
\begin{equation}\label{4.27}
\hat F({\hat W}^R_{k,m},\Phi_{k+1,m}, (k+1)h,m\tau)-
\hat F({\hat W}^L_{k,m},\Phi_{k-1,m}, (k-1)h,m\tau)=O(h^2).
\end{equation}
Our algorithm is so designed that \eqref{4.27} holds as it is one of the keys when proving the limit of $\hat W_h$ is the solution to \eqref{Weq2}, \eqref{4.12}.

Next, for any $k$ with $k+m$ even, we define $\hat W_h$ on the rectangle 
\[
\mathcal R_{k,m}=\left\{(x,t): (k-1)h<x<(k+1)h,\,m\tau\leq t<(m+1)\tau \right\}
\]
by 
\begin{equation}\label{4.28}
\hat W_h(x,t)=U(x-kh,t-m\tau)-\Phi_{k,m},\qquad (x,t)\in\mathcal R_{k,m},
\end{equation}
where $U$ is the solution to the Riemann problem
\begin{equation}\label{4.29}
\begin{split}
&\partial_t U(x,t) +\partial_x F(U(x,t),,kh,m\tau)=0,\qquad-\infty<x<\infty,\quad 0<t<\infty\\
&U(x,0)=\begin{cases} W^R_{k-1,m},&-\infty<x<0,\\W^L_{k+1,m},&0<x<\infty. \end{cases}
\end{split}
\end{equation}

We notice that $\hat W_h$ satisfies the equation
\beq\label{4.30}
\partial_t \hat W_h(x,t)+\partial_x \hat F(\hat W_h(x,t),\Phi_{k,m},kh,m\tau)=0
\eeq
on the rectangle $\mathcal R_{k,m}$, together with the initial condition
\beq\label{4.31}
\hat W_h(x,m\tau)=\begin{cases}
\hat W^R_{k-1,m},&(k-1)h<x<kh,
\\
\hat W^L_{k+1,m},&kh<x<(k+1)h
\end{cases}
\end{equation}
along the base of  $\mathcal R_{k,m}$. Thus \eqref{4.30} is the homogeneous system of conservation laws resulting from \eqref{Weq2} by dropping the source term $\hat G$ and freezing $\Phi$, $x$, $t$ in $\hat F$ at its value $\Phi_{k,m}$ at $(kh,m\tau)$, $kh$ and $m\tau$. On the other hand, the initial data in \eqref{4.31} account for the effect of the source and also of the inhomogeneity due to $\Phi$, $x$ and $t$ in $\hat F$ through \eqref{4.24} and \eqref{4.26}. That is the main idea of operator splitting by Dafermos and Hsiao \cite{DH82} and Dafermos \cite{D06}.

\subsection{Main Estimates}
We now carry out the main estimates that lead to bounds on the total variation of the approximate solution $\hat W_h$.

Recall the structure of self-similar solutions of the strictly hyperbolic system of conservation laws \eqref{4.29}. Its characteristic families are genuinely nonlinear, see Section 4.1. Then any outgoing admissible solution may be visualized as a centered wave fan containing two elementary waves, compressible shocks or centered rarefaction waves, one for each characteristic family. Each elementary wave has a signed amplitude, negative for shocks and positive for rarefactions. Its strength is measured by the absolute value of its amplitude. 

The amplitudes of the two elementary waves composing a wave fan are grouped together in the amplitude vector, which will be denoted by one of the letters $\alpha,\beta, \gamma, \epsilon,\zeta, \xi$. The two end-states, left and right, of a wave fan uniquely determine the amplitude vector via the solution of the corresponding Riemann problem. Conversely, either one of the end-states together with the amplitude vector determine uniquely the opposite end-state.

Specifically, there are $C^2$ functions $P(U^L,\gamma)$, $Q(U^R,\gamma)$ and $\Omega(U^L,U^R)$, determined on some neighborhood of the origin in $\mathbb R^2\times\mathbb R^2$ and taking values in $\mathbb R^2$, such that
\begin{equation}\label{4.32}
U^R=P(U^L,\gamma),\qquad U^L=Q(U^R,\gamma),
\qquad \gamma=\Omega(U^L,U^R),
\end{equation}
for any wave fan with left state $U^L$, right state $U^R$ and amplitude vector $\gamma$. These functions and their partial derivatives satisfy the following for any state $\bar U$:
\begin{eqnarray}
&&P(\bar U,0)=\bar U,\quad Q(\bar U,0)=\bar U,\quad \Omega(\bar U,\bar U)=0,\label{4.33}\\
&&P_{U^L}(\bar U,0)=I, \quad  Q_{U^R}(\bar U,0)=I,\label{4.34}\\
&&P_{\gamma}(\bar U,0)=R(\bar U), \quad Q_{\gamma}(\bar U,0)=-R(\bar U),\label{4.35}\\
&&\Omega_{U^L}(\bar U,\bar U)=-R^{-1}(\bar U),\quad \Omega_{U^R}(\bar U,\bar U)=R^{-1}(\bar U).\label{4.36}
\end{eqnarray}
Here $I$ is the $2\times2$ identity, and
$$
R(\bar U)=\begin{pmatrix}1&1\\ \lambda_-(\bar U,kh,m\tau)&\lambda_+(\bar U,kh,m\tau)\end{pmatrix},
$$
see \eqref{4.15}.

For convenience we write 
\begin{equation}\label{4.37}
H(U,Z)\equiv\Omega(U,U+Z).
\end{equation}
Then from \eqref{4.33} and \eqref{4.36} we have
\begin{equation}\label{4.38}
H(\bar U,0)=0,\quad H_U(\bar U,0)=0,\quad H_Z(\bar U,0)=R^{-1}(\bar U).
\end{equation}

To estimate the total variation of $\hat W_h$ we note that $TV\hat W_h(\cdot,t)$ is constant on each time interval $(m\tau,(m+1)\tau)$:
\begin{equation}\label{4.39}
TV\hat W_h(\cdot,t)=K_m+L_m, \qquad m\tau<t<(m+1)\tau,
\end{equation}
where $K_m$ is the portion of total vaiation contributed by the jump discontinuities of $\hat W_h(\cdot,t)$ at points $x=kh$, with $k+m$ odd, while $L_m$ is the portion induced by the elementary waves emanating from the mesh-points $(kh,m\tau)$, for $k+m$ even.

Now we make the ansatz
\begin{equation}\label{4.40}
L_m\le\varepsilon_2, \qquad m=0,1,\dots,m^*,
\end{equation}
where $\varepsilon_2$ is a small positive number to be fixed during our analysis. In the following the generic constant $C$ is independent to $\varepsilon_2$ as well.

First we estimate the jump of $\hat W_h(\cdot,t)$ across $x=kh$, with $k+m$ odd. From \eqref{4.31}, \eqref{4.26} and \eqref{4.24}, it is
\begin{equation}\label{4.41}
\hat W^R_{k,m}-\hat W^L_{k,m}=(S_{k,m}-I)(\Phi_{k+1,m}-\Phi_{k-1,m})-2h\tilde S_{k,m}.
\end{equation}
From \eqref{4.21} and \eqref{4.20},
\begin{equation}\label{4.42}
S_{k,m}-I=O(1)|W_{k,m}-\Phi_{k,m}|=O(1)|\hat W_h(y_{k,m},m\tau-)|.
\end{equation}
From \eqref{4.25},
\begin{equation}\label{4.43}
\tilde S_{k,m}=O(1)|\hat W_{k,m}\theta_x(kh,m\tau)|=O(1)|\hat W_h(y_{k,m},m\tau-)||\theta_x(kh,m\tau)|.
\end{equation}

Summing up \eqref{4.41} for $k$ with $k+m$ odd and substituting \eqref{4.42} and \eqref{4.43} into it, we arrive at
\begin{equation*}
\begin{split}
K_m&=\sum|\hat W^R_{k,m}-\hat W^L_{k,m}|\le C\sum|\hat W_h(y_{k,m},m\tau-)|[|\Phi_{k+1,m}-\Phi_{k-1,m}|+2h|\theta_x(kh,m\tau)|]\\
&\le C\,TV\hat W_h(\cdot,m\tau-)\sum[|\Phi_{k+1,m}-\Phi_{k-1,m}|+2h|\theta_x(kh,m\tau)|].
\end{split}
\end{equation*}
Applying \eqref{4.18}, \eqref{4.5}, \eqref{4.3}, \eqref{l1w}, \eqref{thdef} and \eqref{Mes} and by iteration, we have
\begin{equation}\label{4.44}
K_m\le C(K_{m-1}+L_{m-1})\sigma(m\tau+1)^{-\frac14}\le C\varepsilon_2\sigma(m\tau+1)^{-\frac14}.
\end{equation}
Here in the iteration we have taken $\delta_0$ and $\sigma_0$ in Theorem \ref{main_thm} small so they are bounded by $\varepsilon_2$.

With \eqref{4.40} and \eqref{4.44}, it is straightforward to verify the following,
\begin{equation}\label{4.45}
|W_{k,m}|,\,|\hat W_{k,m}|,\,|W_{k,m}^L|,\,|W_{k,m}^R|,\,|\hat W_{k,m}^L|,\,|\hat W_{k,m}^R|\le C\varepsilon_2.
\end{equation}
We also update \eqref{4.42} and \eqref{4.43} to
\begin{eqnarray}\label{4.46a}
&&S_{k,m}-I=O(1)\varepsilon_2,\\
&&\tilde S_{k,m}=O(1)\varepsilon_2|\theta_x(kh,m\tau)|.\label{4.46}
\end{eqnarray}

Our next step is to estimate $L_m$. It is measured as the sum of the strengths of elementary waves from the mesh-points $(kh,m\tau)$, $k+m$ even, and cross the $t$-time line, for $t\in(m\tau,(m+1)\tau)$. To estimate $L_m$ inductively, we follow the standard approach in the random choice method to compare the strengths of incoming waves and outgoing waves in a typical diamond-shaped domain. Precisely, for $k+m$ even, we connect the random mesh-points $(y_{k,m-1},(m-1)\tau)$, $(y_{k-1,m},m\tau)$, $(y_{k,m+1},(m+1)\tau)$ and $(y_{k+1,m},m\tau)$ to form a region $\mathcal  D_{k,m}$, see the diamond-shaped domain formed with dashed lines in Figure 1.

\begin{figure}
\begin{center}
\includegraphics[width=0.6\textwidth,height=0.3\textheight]{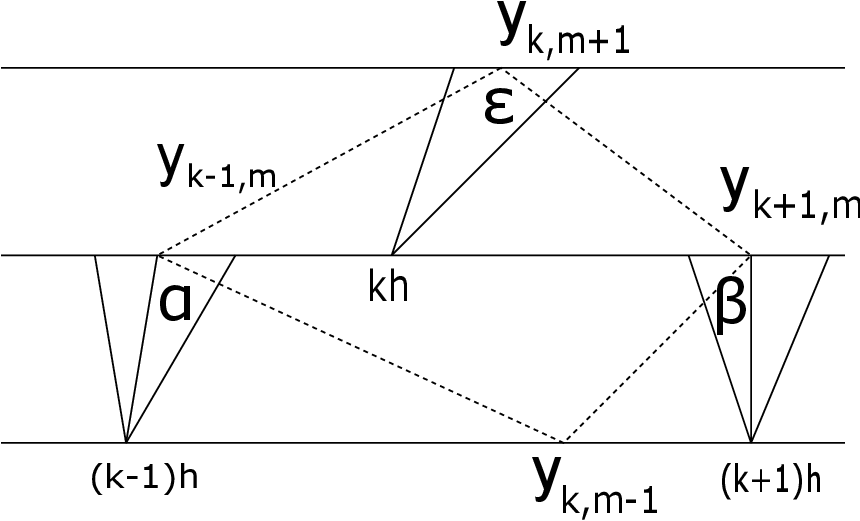}
\caption{The random choice algorithm: $m+k$ is even. $\alpha$, $\beta$ and $\epsilon$ are corresponding wave strengths.}
\end{center}
\end{figure}

Let $\alpha$ be the part of the wave fan emanating from the mesh-point $((k-1)h,(m-1)\tau)$ that enters $\mathcal  D_{k,m}$ through its southwestern side, and  $\beta$ be the part of the wave fan emanating from the mesh-point $((k+1)h,(m-1)\tau)$ that enters $\mathcal  D_{k,m}$ through its southeastern side.
On the other hand, the outgoing waves consists of the full wave fan, with amplitude vector $\epsilon$, emanating from the mesh-point $(kh,m\tau)$ and exiting through the northwestern and/or northeastern side of $\mathcal  D_{k,m}$.
	
By \eqref{4.20}, \eqref{4.28}, \eqref{4.29} and \eqref{4.32},  we write
\begin{eqnarray}
&&W_{k-1,m}=Q(W^L_{k,m-1},\alpha)+\Phi_{k-1,m}-\Phi_{k-1,m-1},\label{4.47}\\
&&W_{k+1,m}=P(W^R_{k,m-1},\beta)+\Phi_{k+1,m}-\Phi_{k+1,m-1},\label{4.48}\\
&&\epsilon=\Omega(W^R_{k-1,m},W^L_{k+1,m})=H(W^R_{k-1,m},W^L_{k,+1m}-W^R_{k-1,m}).\label{4.49}
\end{eqnarray}
	
We now define
\begin{eqnarray}
&&\tilde W_{k-1,m}=	Q(W_{k,m-1},\alpha),\qquad\tilde W_{k+1,m}=P(W_{k,m-1},\beta),\label{4.50}\\
&&\tilde\epsilon=\Omega(\tilde W_{k-1,m},\tilde W_{k+1,m})=H(\tilde W_{k-1,m},\tilde W_{k+1,m}-\tilde W_{k-1,m}).\label{4.51}
\end{eqnarray}
By Glimm's celebrated wave interaction estimates \cite{Glimm},
\begin{equation}\label{4.52}
\tilde\epsilon=\alpha+\beta+O(X_{k,m}),
\end{equation}
where
$$
X_{k,m}=\sum|\alpha_i||\beta_j|,
$$
with the summation running over all pairs $(i,j)$ of approaching elementary waves of the two incoming wave fans.

Next we estimate $\epsilon-\tilde\epsilon$ following the approach in \cite{D06}. We use \eqref{4.49} and \eqref{4.51} to write
\begin{equation}\label{4.53}
\epsilon-\tilde\epsilon=\bar H_U(W^R_{k-1,m}-\tilde W_{k-1,m})+\bar H_Z(W^L_{k+1,m}-W^R_{k-1,m}-\tilde W_{k+1,m}+\tilde W_{k-1,m}),
\end{equation}
where
\begin{eqnarray}
&&\bar H_U=\int_0^1H_U(sW^R_{k-1,m}+(1-s)\tilde W_{k-1,m}, \tilde W_{k+1,m}-\tilde W_{k-1,m})\,ds,\label{4.54}\\
&&\bar H_Z=\int_0^1H_Z(W^R_{k-1,m},s(W^L_{k+1,m}-W^R_{k-1,m})+(1-s)(\tilde W_{k+1,m}- \tilde W_{k-1,m}))\,ds.\label{4.55}
\end{eqnarray}
We further use \eqref{4.24}, \eqref{4.50}, \eqref{4.47}, \eqref{4.48}
to write
\begin{eqnarray}
&W^R_{k-1,m}-\tilde W_{k-1,m}=&\bar Q_{U^R}[S_{k,m-1}(\Phi_{k-1,m-1}-\Phi_{k,m-1})+h\tilde S_{k,m-1}-\tau G_{k,m-1}]\label{4.56}\\
&&+\Phi_{k-1,m}-\Phi_{k-1,m-1}
+S_{k-1,m}(\Phi_{k,m}-\Phi_{k-1,m})\nonumber\\
&&-h\tilde S_{k-1,m}-\tau G_{k-1,m},\nonumber\\
&W^L_{k+1,m}-\tilde W_{k+1,m}=&\bar P_{U^L}[S_{k,m-1}(\Phi_{k+1,m-1}-\Phi_{k,m-1})-h\tilde S_{k,m-1}-\tau G_{k,m-1}]\label{4.57}\\
&&+\Phi_{k+1,m}-\Phi_{k+1,m-1}
+S_{k+1,m}(\Phi_{k,m}-\Phi_{k+1,m})\nonumber\\
&&+h\tilde S_{k+1,m}-\tau G_{k+1,m},\nonumber
\end{eqnarray}
where
\begin{equation}\label{4.58}
\begin{split}
&\bar Q_{U^R}=\int_0^1Q_{U^R}(sW^L_{k,m-1}+(1-s)W_{k,m-1},\alpha)\, ds,\\
&\bar P_{U^L}=\int_0^1P_{U^L}(sW^R_{k,m-1}+(1-s)W_{k,m-1},\beta)\, ds.
\end{split}
\end{equation}
We note that those terms with $\tilde S_{k,m-1}$ and $\tilde S_{k\pm1,m}$ are new terms when comparing with \cite{D06}.

From \eqref{4.54}, \eqref{4.38}, \eqref{4.50} and \eqref{4.33} one has
$$
|\bar H_U|\le C|\tilde W_{k+1,m}-\tilde W_{k-1,m}|\le C(|\alpha|+|\beta|).
$$
Taking account of \eqref{4.46} in \eqref{4.56}, we have the similar estimate as in \cite{D06} for the first term on the right-hand side of \eqref{4.53}:
\begin{equation}\label{4.59}
|\bar H_U(W^R_{k-1,m}-\tilde W_{k-1,m})|\le C(\varepsilon_1+\varepsilon_2)(|\alpha|+|\beta|)h.
\end{equation}
Here we have applied \eqref{4.34}, \eqref{4.46a}, \eqref{4.2}, \eqref{4.3}, \eqref{4.5}, \eqref{4.23}, \eqref{Weq4}, \eqref{4.45}, \eqref{l1w} and \eqref{4.4}, and taking $\sigma_0\le \varepsilon_2$.

The estimate of the second term on the right-hand side of \eqref{4.53} is much more complicated. First, with \eqref{4.56} and \eqref{4.57} we write
\begin{equation}\label{4.61}
\begin{split}
&W^L_{k+1,m}-W^R_{k-1,m}-\tilde W_{k+1,m}+\tilde W_{k-1,m}\\
=&-\tau(G_{k+1,m}-G_{k-1,m})+h(\tilde S_{k+1,m}+\tilde S_{k-1,m}-2\tilde S_{k,m-1})\\
&+(\bar P_{U^L}-I)[S_{k,m-1}(\Phi_{k+1,m-1}-\Phi_{k,m-1})-h\tilde S_{k,m-1}-\tau G_{k,m-1}]\\
&-(\bar Q_{U^R}-I)[S_{k,m-1}(\Phi_{k-1,m-1}-\Phi_{k,m-1})-h\tilde S_{k,m-1}-\tau G_{k,m-1}]\\
&-(S_{k,m-1}-I)(\Phi_{k+1,m}-\Phi_{k+1,m-1}+\Phi_{k-1,m-1}-\Phi_{k-1,m})\\
&+(S_{k+1,m}-S_{k,m-1})(\Phi_{k,m}-\Phi_{k+1,m})+(S_{k,m-1}-S_{k-1,m})(\Phi_{k,m}-\Phi_{k-1,m}).
\end{split}
\end{equation}
We focus on the first two terms on the right-hand side since the former is the leading term and related to our specific $\hat G_{\hat W}$ in \eqref{4.16} while the latter is a new contribution.

For the first term on the right-hand side of \eqref{4.61}, by \eqref{4.23} we further write
\begin{equation}\label{4.62}
G_{k+1,m}-G_{k-1,m}=\overline{\hat G}_{\hat W}(\hat W_{k+1,m}-\hat W_{k-1,m})+\overline{\hat G}_{\Phi}(\Phi_{k+1,m}-\Phi_{k-1,m})+2h\overline{\hat G}_{x},
\end{equation}
where
\begin{equation}\label{4.63}
\begin{split}
&\overline{\hat G}_{\hat W}=\int_0^1\hat G_{\hat W}(s\hat W_{k+1,m}+(1-s)\hat W_{k-1,m},\Phi_{k+1,m},(k+1)h,m\tau)\,ds,\\
&\overline{\hat G}_{\Phi}=\int_0^1\hat G_{\Phi}(\hat W_{k-1,m},s\Phi_{k+1,m}+(1-s)\Phi_{k-1,m},(k+1)h,m\tau)\,ds,\\
&\overline{\hat G}_{x}=\int_0^1\hat G_x(\hat W_{k-1,m},\Phi_{k-1,m},(k-1)h+2sh,m\tau)\,ds.
\end{split}
\end{equation}
Here, by \eqref{4.20},  \eqref{4.47}, \eqref{4.48}, \eqref{4.33} and \eqref{4.24} we have
\begin{equation}\label{4.64}
\hat W_{k+1,m}-\hat W_{k-1,m}=\bar P_\gamma\beta-\bar Q_\gamma\alpha+(S_{k,m-1}-I)(\Phi_{k+1,m-1}-\Phi_{k-1,m-1})-2h\tilde S_{k,m-1},
\end{equation}
where
\begin{equation}\label{4.65}
\bar P_\gamma=\int_0^1P_\gamma(W^R_{k,m-1},s\beta)\,ds,\qquad \bar Q_\gamma=\int_0^1Q_\gamma(W^L_{k,m-1},s\alpha)\,ds.
\end{equation}
Thus, substituting \eqref{4.64} into \eqref{4.62} and noting $\overline{\hat G}_{\hat W}$ and $\overline{\hat G}_{\Phi}$ are $O(1)$ due to \eqref{4.17}, together with \eqref{4.46a},  
the first term on the right-hand side of \eqref{4.61} becomes
\begin{equation}\label{4.66}
\begin{split}
-\tau(G_{k+1,m}-G_{k-1,m})=&\,\tau\overline{\hat G}_{\hat W}\bar Q_\gamma\alpha-\tau\overline{\hat G}_{\hat W}\bar P_\gamma\beta+O(1)\varepsilon_2\tau|\Phi_{k+1,m-1}-\Phi_{k-1,m-1}|\\
&+O(1)\tau h|\tilde S_{k,m-1}|+O(1)\tau|\Phi_{k+1,m}-\Phi_{k-1,m}|+O(1)\tau h|\overline{\hat G}_{x}|.
\end{split}
\end{equation}

We now consider the second term on the right-hand side of \eqref{4.61}. We write
\begin{equation}\label{4.67}
|\tilde S_{k+1,m}+\tilde S_{k-1,m}-2\tilde S_{k,m-1}|\le|\tilde S_{k+1,m}-\tilde S_{k,m-1}|+|\tilde S_{k-1,m}-\tilde S_{k,m-1}|.
\end{equation}
From \eqref{4.25} one has
\begin{equation}\label{4.68}
\begin{split}
|\tilde S_{k+1,m}-\tilde S_{k,m-1}|\le& C\varepsilon_2|W_{k+1,m}-W_{k,m-1}||\theta_x((k+1)h,m\tau)|\\
&+C\varepsilon_2||\theta_x((k+1)h,m\tau)-\theta_x(kh,(m-1)\tau)|\\
&+C(|\hat W_{k+1,m}-\hat W_{k,m-1}||\theta_x(kh,(m-1)\tau)|.
\end{split}
\end{equation}
Applying \eqref{4.48} and \eqref{4.33}  we have
\begin{equation}\label{4.69}
\begin{split}
|W_{k+1,m}-W_{k,m-1}|\le& |P(W^R_{k,m-1},\beta)-W_{k,m-1}|+|\Phi_{k+1,m}-\Phi_{k+1,m-1}|\\
\le&C|\beta|+|W^R_{k,m-1}-W_{k,m-1}|+|\Phi_{k+1,m}-\Phi_{k+1,m-1}|.
\end{split}
\end{equation}
Applying \eqref{4.18} and \eqref{4.2}-\eqref{4.5} one can verify
\begin{equation}\label{4.70}
\begin{split}
&|\Phi_{k,m}-\Phi_{k,m-1}|\le C(\varepsilon_1+\sigma)\tau\le C(\varepsilon_1+\varepsilon_2)\tau,\\
&|\Phi_{k+1,m}-\Phi_{k,m}|\le C\varepsilon_1h,
\end{split}
\end{equation}
taking $\sigma_0\le\varepsilon_2$. Together with \eqref{4.24},  \eqref{4.46a}, \eqref{4.46},  \eqref{4.23}, \eqref{Weq4} and \eqref{4.45},  from      \eqref{4.69}  we deduce
\begin{equation}\label{4.71}
|W_{k+1,m}-W_{k,m-1}|\le C|\beta|+C(\varepsilon_1+\varepsilon_2)\tau.
\end{equation}
Similarly, one can find
\begin{equation}\label{4.72}
|\hat W_{k+1,m}-\hat W_{k,m-1}|\le C|\beta|+C(\varepsilon_1+\varepsilon_2)\tau.
\end{equation}
Also, from the formulation of $\theta$, \eqref{thdef} and \eqref{Mes}, one has
\begin{equation}\label{4.73}
|\theta_x((k+1)h,m\tau)-\theta_x(kh,(m-1)\tau)|\le C\sigma\tau(m\tau+1)^{-\frac32}e^{-\frac{(kh)^2}{8(m\tau+1)}}.
\end{equation}
Substituting \eqref{4.71}-\eqref{4.73} into \eqref{4.68}, we arrive at
\begin{equation}\label{4.74}
\begin{split}
|\tilde S_{k+1,m}-\tilde S_{k,m-1}|\le& C\sigma|\beta|+C(\varepsilon_1+\varepsilon_2)\tau[|\theta_x((k+1)h,m\tau)|+|\theta_x(kh,(m-1)\tau)|]\\
&+C\varepsilon_2\sigma\tau(m\tau+1)^{-\frac32}e^{-\frac{(kh)^2}{8(m\tau+1)}}.
\end{split}
\end{equation}

Similar to \eqref{4.74} we also have
\begin{equation}\label{4.75}
\begin{split}
|\tilde S_{k-1,m}-\tilde S_{k,m-1}|\le& C\sigma|\alpha|+C(\varepsilon_1+\varepsilon_2)\tau[|\theta_x((k-1)h,m\tau)|+|\theta_x(kh,(m-1)\tau)|]\\
&+C\varepsilon_2\sigma\tau(m\tau+1)^{-\frac32}e^{-\frac{(kh)^2}{8(m\tau+1)}}.
\end{split}
\end{equation}
Substituting \eqref{4.74} and \eqref{4.75} into \eqref{4.67}, we have the bound of the second term on the right-hand side of \eqref{4.61} as
\begin{equation}\label{4.76}
\begin{split}
h|\tilde S_{k+1,m}+\tilde S_{k-1,m}-2\tilde S_{k,m-1}|\le&C\sigma(|\alpha|+|\beta|)h+C(\varepsilon_1+\varepsilon_2)\tau h[|\theta_x((k+1)h,m\tau)|\\
&+|\theta_x((k-1)h,m\tau)|+|\theta_x(kh,(m-1)\tau)|]\\
&+C\varepsilon_2\sigma\tau h(m\tau+1)^{-\frac32}e^{-\frac{(kh)^2}{8(m\tau+1)}}.
\end{split}
\end{equation}

The other terms on the right-hand side of \eqref{4.61} are parallel to those in \cite{D06}. With slight modifications due to the terms containing $\tilde S_{k,m-1}$, they can be estimated in a similar way. This gives us the following,
\begin{equation}\label{4.77}
\begin{split}
&\left|(\bar P_{U^L}-I)[S_{k,m-1}(\Phi_{k+1,m-1}-\Phi_{k,m-1})-h\tilde S_{k,m-1}-\tau G_{k,m-1}]\right|
\le C|\beta|(\varepsilon_1+\varepsilon_2)\tau ,\\
&\left|(\bar Q_{U^R}-I)[S_{k,m-1}(\Phi_{k-1,m-1}-\Phi_{k,m-1})-h\tilde S_{k,m-1}-\tau G_{k,m-1}]\right|
\le C|\alpha|(\varepsilon_1+\varepsilon_2)\tau, \\
&\left|(S_{k+1,m}-S_{k,m-1})(\Phi_{k,m}-\Phi_{k+1,m})\right|\le C[|\beta|+(\varepsilon_1+\varepsilon_2)\tau]|\Phi_{k,m}-\Phi_{k+1,m}|, \\
&\left|(S_{k,m-1}-S_{k-1,m})(\Phi_{k,m}-\Phi_{k-1,m})\right|
\le C[|\alpha|+(\varepsilon_1+\varepsilon_2)\tau ]|\Phi_{k,m}-\Phi_{k-1,m}|,
\end{split}
\end{equation}
applying \eqref{4.58}, \eqref{4.34}, \eqref{4.45}-\eqref{4.46}, \eqref{4.23}, \eqref{Weq4}, \eqref{4.21}, \eqref{4.70} and \eqref{4.71}.

Now we substitute \eqref{4.66}, \eqref{4.76} and \eqref{4.77} into the right-hand side of \eqref{4.61}. Together with \eqref{4.59}, we give an estimate for $\epsilon-\tilde \epsilon$ in \eqref{4.53} as follows,
\begin{equation}\label{4.78}
\begin{split}
\epsilon-\tilde \epsilon=&\,
\tau\bar H_Z\overline{\hat G}_{\hat W}\bar Q_\gamma\alpha-\tau\bar H_Z\overline{\hat G}_{\hat W}\bar P_\gamma\beta+O(1)(\varepsilon_1+\varepsilon_2)\tau(|\alpha|+|\beta|)\\
&+O(1)\varepsilon_2|\Phi_{k+1,m}-\Phi_{k+1,m-1}+\Phi_{k-1,m-1}-\Phi_{k-1,m}|\\
&+O(1)\tau [h|\tilde S_{k,m-1}|+|\Phi_{k+1,m}-\Phi_{k,m}|+|\Phi_{k,m}-\Phi_{k-1,m}|]\\
&+O(1)(\varepsilon_1+\varepsilon_2)\tau h[|\theta_x((k+1)h,m\tau)|
+|\theta_x((k-1)h,m\tau)|+|\theta_x(kh,(m-1)\tau)|]\\
&+O(1)\tau h\max_{[(k-1)h,(k+1)h]}(|\theta_x|+|\theta_{xx}|+|\theta_{xxx}|)(\cdot,m\tau)\\
&+O(1)\varepsilon_2\sigma\tau h(m\tau+1)^{-\frac32}e^{-\frac{(kh)^2}{8(m\tau+1)}},
\end{split}
\end{equation}
where we have applied \eqref{4.55}, \eqref{4.38}, \eqref{4.46a}, \eqref{4.17} and \eqref{4.63}.

Combining \eqref{4.78} with \eqref{4.52} gives us
\begin{equation}\label{4.79}
\begin{split}
\epsilon=&\,
(I+\tau\bar H_Z\overline{\hat G}_{\hat W}\bar Q_\gamma)\alpha+(I-\tau\bar H_Z\overline{\hat G}_{\hat W}\bar P_\gamma)\beta+O(1)(\varepsilon_1+\varepsilon_2)\tau(|\alpha|+|\beta|)\\
&+O(1)\varepsilon_2|\Phi_{k+1,m}-\Phi_{k+1,m-1}+\Phi_{k-1,m-1}-\Phi_{k-1,m}|\\
&+O(1)\tau [h|\tilde S_{k,m-1}|+|\Phi_{k+1,m}-\Phi_{k,m}|+|\Phi_{k,m}-\Phi_{k-1,m}|]\\
&+O(1)(\varepsilon_1+\varepsilon_2)\tau h[|\theta_x((k+1)h,m\tau)|
+|\theta_x((k-1)h,m\tau)|+|\theta_x(kh,(m-1)\tau)|]\\
&+O(1)\tau h\max_{[(k-1)h,(k+1)h]}(|\theta_x|+|\theta_{xx}|+|\theta_{xxx}|)(\cdot,m\tau)\\
&+O(1)\varepsilon_2\sigma\tau h(m\tau+1)^{-\frac32}e^{-\frac{(kh)^2}{8(m\tau+1)}}+O(X_{k,m}).
\end{split}
\end{equation}
By the specific form of $\hat G_{\hat W}$ in \eqref{4.16}, together with \eqref{4.55}, \eqref{4.63}, \eqref{4.65}, \eqref{4.38}, \eqref{4.35}, \eqref{4.45}, \eqref{4.50} and \eqref{4.33}, we have
\begin{equation}\label{4.80}
\bar H_Z\overline{\hat G}_{\hat W}\bar Q_\gamma=-\frac12I+O(1)\varepsilon_2,\qquad\bar H_Z\overline{\hat G}_{\hat W}\bar P_\gamma=\frac12I+O(1)\varepsilon_2,
\end{equation}
taking $\sigma_0\le\varepsilon_2$.
Therefore, taking $\varepsilon_2$ small, there exists a positive constant $\nu$ such that
\begin{equation}\label{4.81}
|I+\tau\bar H_Z\overline{\hat G}_{\hat W}\bar Q_\gamma|\le1-4\nu\tau,\qquad|I-\tau\bar H_Z\overline{\hat G}_{\hat W}\bar P_\gamma|\le1-4\nu\tau.
\end{equation}
Substituting \eqref{4.81} into \eqref{4.79} and taking $\varepsilon_1$ and $\varepsilon_2$ small, we arrive at
\begin{equation}\label{4.82}
\begin{split}
|\epsilon|\le&\,
(1-3\nu\tau)(|\alpha|+|\beta|)
+C\varepsilon_2|\Phi_{k+1,m}-\Phi_{k+1,m-1}+\Phi_{k-1,m-1}-\Phi_{k-1,m}|\\
&+C\tau [h|\tilde S_{k,m-1}|+|\Phi_{k+1,m}-\Phi_{k,m}|+|\Phi_{k,m}-\Phi_{k-1,m}|]\\
&+C(\varepsilon_1+\varepsilon_2)\tau h[|\theta_x((k+1)h,m\tau)|
+|\theta_x((k-1)h,m\tau)|+|\theta_x(kh,(m-1)\tau)|]\\
&+C\tau h\max_{[(k-1)h,(k+1)h]}(|\theta_x|+|\theta_{xx}|+|\theta_{xxx}|)(\cdot,m\tau)\\
&+C\varepsilon_2\sigma\tau h(m\tau+1)^{-\frac32}e^{-\frac{(kh)^2}{8(m\tau+1)}}+O(X_{k,m}).
\end{split}
\end{equation}

We can obtain an estimate similar to \eqref{4.79} for $\epsilon-\alpha-\beta$. Then we also have
\begin{equation}\label{4.83}
\begin{split}
|\epsilon-\alpha-\beta|\le&\,
C\tau(|\alpha|+|\beta|)
+C\varepsilon_2|\Phi_{k+1,m}-\Phi_{k+1,m-1}+\Phi_{k-1,m-1}-\Phi_{k-1,m}|\\
&+C\tau [h|\tilde S_{k,m-1}|+|\Phi_{k+1,m}-\Phi_{k,m}|+|\Phi_{k,m}-\Phi_{k-1,m}|]\\
&+C(\varepsilon_1+\varepsilon_2)\tau h[|\theta_x((k+1)h,m\tau)|
+|\theta_x((k-1)h,m\tau)|+|\theta_x(kh,(m-1)\tau)|]\\
&+C\tau h\max_{[(k-1)h,(k+1)h]}(|\theta_x|+|\theta_{xx}|+|\theta_{xxx}|)(\cdot,m\tau)\\
&+C\varepsilon_2\sigma\tau h(m\tau+1)^{-\frac32}e^{-\frac{(kh)^2}{8(m\tau+1)}}+O(X_{k,m}).
\end{split}
\end{equation}

Holding $m$ fixed, we sum \eqref{4.82} and \eqref{4.83} over all $k$ with $k+m$ even. The sum of the $|\epsilon|$ terms is $L_m$ while the sum of the $(|\alpha|+|\beta|)$ terms is $L_{m-1}$. We denote the sum of the $|\epsilon-\alpha-\beta|$ terms by $\Delta_m$. Note that by \eqref{4.46}, \eqref{4.18}, \eqref{4.3}, \eqref{4.5}, \eqref{l1w} and the specific form of $\theta$ given in \eqref{thdef} and \eqref{Mes} we have
\begin{equation}\label{4.84}
\begin{split}
&\sum_kh|\tilde S_{k,m-1}|\le C\varepsilon_2\sigma(m\tau+1)^{-\frac12},\\
&\sum_k(|\Phi_{k+1,m}-\Phi_{k,m}|+|\Phi_{k,m}-\Phi_{k-1,m}|)\le C\sigma(m\tau+1)^{-\frac14},\\
&\sum_kh
[|\theta_x((k+1)h,m\tau)|
+|\theta_x((k-1)h,m\tau)|+|\theta_x(kh,(m-1)\tau)|]
\le C\sigma(m\tau+1)^{-\frac12},\\
&\sum_kh\max_{[(k-1)h,(k+1)h]}(|\theta_x|+|\theta_{xx}|+|\theta_{xxx}|)(\cdot,m\tau)\le C\sigma(m\tau+1)^{-\frac12},\\
&\sum_kh(m\tau+1)^{-\frac32}e^{-\frac{(kh)^2}{8(m\tau+1)}}\le C(m\tau+1)^{-1}.
\end{split}
\end{equation}
With \eqref{4.4} we also have
\begin{equation}\label{4.85}
\begin{split}
&\sum_k|\Phi_{k+1,m}-\Phi_{k+1,m-1}+\Phi_{k-1,m-1}-\Phi_{k-1,m}|\\
\le&\,\frac12\int_{(m-1)\tau}^{m\tau}\sum_k|w_2((k+1)h,t)-w_2((k-1)h,t)|\,dt\\
&+\frac12\int_{(m-1)\tau}^{m\tau}\sum_k|\theta_x((k+1)h,t)-\theta_x((k-1)h,t)|\,dt\\
\le&\,\frac12J_m+C\sigma\tau(m\tau+1)^{-1},
\end{split}
\end{equation}
where
\begin{equation}\label{4.86}
J_m=\int_{(m-1)\tau}^{m\tau}TVW(\cdot,t)\,dt\le\varepsilon_1\tau.
\end{equation}

Thus, summing up \eqref{4.82} gives us
\begin{equation}\label{4.87}
L_m\le(1-3\nu\tau)L_{m-1}+C\sigma\tau(m\tau+1)^{-\frac14}+C\varepsilon_2J_m+C\sum_{k,\,k+m\,even}X_{k,m},
\end{equation}
and summing up \eqref{4.83} gives us
\begin{equation}\label{4.88}
\Delta_m\le C\tau L_{m-1}+C\varepsilon_2\tau+C\sum_{k,\,k+m\,even}X_{k,m}\le C\varepsilon_2\tau+C\sum_{k,\,k+m\,even}X_{k,m},
\end{equation}
noting \eqref{4.86} and taking $\sigma_0\le\varepsilon_2$.

To control the wave interaction terms $X_{k,m}$ we follow the standard procedure as follows (see \cite{D06} for details). Let
$$
M_m=\sum|\zeta_i||\xi_j|,
$$
where the summation is for all pairs of approaching elementary waves emanating from mesh points along the $m\tau$-time line. It satisfies
\begin{eqnarray} 
&M_m&\le L^2_m\le \varepsilon_2 L_m,\label{4.89}\\
&M_m&\le M_{m-1}+(L_m+L_{m-1})\Delta_m-\sum_{k,\,k+m\,even}X_{k,m}\label{4.90}\\
&&\le  M_{m-1}+(L_m+L_{m-1})C\varepsilon_2\tau+(C\varepsilon_2-1) \sum_{k,\,k+m\,even}X_{k,m}.\nonumber
\end{eqnarray} 

Next, we introduce the Glimm functional
\begin{equation}\label{4.91}
N_m=L_m+\kappa M_m,
\end{equation}
where $\kappa$ is a positive constant. Substituting \eqref{4.87} and \eqref{4.90} into \eqref{4.91} and applying \eqref{4.89}, after simplification we choose a large $\kappa$ and make $\varepsilon_2$ small to arrive at
\begin{equation}\label{4.92}
N_m\le (1-2\nu\tau)N_{m-1}+C\sigma\tau(m\tau+1)^{-\frac14}+C\varepsilon_2J_m.
\end{equation}
After iteration and noting $N_0\le C(\delta+\sigma)$, we have
\begin{equation}\label{4.93}
N_m
\le C(\delta+\sigma)(1-2\nu\tau)^m+C\sigma\tau\sum_{l=1}^m(1-2\nu\tau)^{m-l}(l\tau+1)^{-\frac14}+C\varepsilon_2\sum_{l=1}^m(1-2\nu\tau)^{m-l}J_l.
\end{equation}

From \eqref{4.86} we have
\begin{equation}\label{4.94}
\begin{split}
L_m\le N_m\le& C(\delta+\sigma)+C\sigma\tau\sum_{l=1}^m(1-2\nu\tau)^{m-l}+C\varepsilon_2\sum_{l=1}^m(1-2\nu\tau)^{m-l}\varepsilon_1\tau\\
&\le C(\delta+\sigma)+C\varepsilon_1\varepsilon_2,
\end{split}
\end{equation}
taking $\sigma_0\le \varepsilon_2$. Noting that by choosing $\delta_0,\sigma_0$ and $\varepsilon_1$ small, 
the right-hand side of \eqref{4.94} is bounded be $\varepsilon_2$. This justifies the ansatz in  \eqref{4.40}.

Now from \eqref{4.39}, \eqref{4.44}, \eqref{4.93} and \eqref{4.86} we obtain our key estimate on the total variation of $\hat W_h$ as follows. For $m\tau<t<(m+1)\tau$, $m=0,1,\dots, m^*$,
\begin{eqnarray}
TV\hat W_h(\cdot,t)\le K_m+N_m&\le&C\varepsilon_2\sigma(m\tau+1)^{-\frac14}+C(\delta+\sigma)(1-2\nu\tau)^m\label{4.95}\\
&&+C\sigma\tau\sum_{l=1}^m(1-2\nu\tau)^{m-l}(l\tau+1)^{-\frac14}\nonumber\\
&&+C\varepsilon_2\sum_{l=1}^m(1-2\nu\tau)^{m-l}\int_{(l-1)\tau}^{l\tau}TVW(\cdot,t)\,dt\nonumber\\
&\le &C(\delta+\sigma)+C\varepsilon_1\varepsilon_2.\label{4.96}
\end{eqnarray}

Next, we fix $s$ and $t$ with, say, $(j-1)\tau<s<j\tau\le l\tau<t\le(l+1)\tau$, and some large positive number $R$. With minor adjustments related to $\tilde S_{k,m}$ and the dependence of $G_{k,m}$ on $\theta$ and its derivatives,  one can show that
\begin{equation}\label{4.97}
\int_{-R}^R|\hat W_h(x,t)-\hat W_h(x,s)|\,dx\le CR\varepsilon_2[|t-s|+h].
\end{equation}
Readers are referred to \cite{D06} for details.

\subsection{Prood of Theorem \ref{main_thm}}
We are now ready to prove our main result.
With the estimates \eqref{4.96} and \eqref{4.97}, we apply Helley's Theorem and go through a standard diagonal process to extract a sequence $\{h_l\}$, with $\lim_{l\rightarrow\infty}h_l=0$, and a function $\hat W$ defined on the strip $(-\infty,\infty)\times[0,T)$, such that
\begin{equation}\label{4.98}
\lim_{l\rightarrow\infty}\hat W_{h_l}(x,t)=\hat W(x,t)
\end{equation}
for all $t\in[0,T)$ and almost all $x\in(-\infty,\infty)$. From \eqref{4.95} we further have
\begin{equation}\label{4.99}
\begin{split}
TV\hat W(\cdot,t)\le &C\varepsilon_2\sigma(t+1)^{-\frac14}+C(\delta+\sigma)e^{-2\nu t}
+C\sigma\int_0^te^{-2\nu(t-s)}(s+1)^{-\frac14}\,ds\\
&+C\varepsilon_2\int_0^te^{-2\nu(t-s)}TVW(\cdot,s)\,ds\\
\le &C\sigma(t+1)^{-\frac14}+C\delta e^{-2\nu t}+C\varepsilon_2\int_0^te^{-2\nu(t-s)}TVW(\cdot,s)\,ds
\end{split}
\end{equation}
for $t\in[0,T)$.

Next, one can show that $\hat W$ is a solution of \eqref{Weq2}, \eqref{4.12}, see \cite{D06} for details. Here 
a key estimate from Glimm's random choice method and \eqref{4.27} play essential roles. Then uniqueness of admissible solutions to the Cauchy problem implies that the $\hat W$ in \eqref{4.98} is the same as the $\hat W$ related to $W$ through \eqref{Weq1}, with $\Phi$ defined by \eqref{4.3} and \eqref{4.5}. In particular, we combine \eqref{Weq1}, \eqref{4.99} and \eqref{l1w} to have
\begin{equation}\label{4.100}
\begin{split}
TV W(\cdot,t)&\le TV\hat W(\cdot,t)+TV\Phi(\cdot,t)\\
&\le C\sigma(t+1)^{-\frac14}+C\delta e^{-2\nu t}+C\varepsilon_2\int_0^te^{-2\nu(t-s)}TVW(\cdot,s)\,ds
\end{split}
\end{equation}
for $t\in[0,T)$.

Noting \eqref{4.100} is a Gronwall's-type integral inequality, we take $\varepsilon_2$ small and solve the inequality to obtain
\begin{equation}\label{4.101}
TV W(\cdot,t)\le C\sigma(t+1)^{-\frac14}+C\delta e^{-\nu t},\quad t\in[0,T).
\end{equation}
By taking $\sigma_0$ and $\delta_0$ small, \eqref{4.101} justifies our assumption \eqref{4.2}. Finally, recalling \eqref{defw}, we further have
\begin{equation}\label{4.102}
TV v(\cdot,t)+TV\tilde u(\cdot,t)\le TV W(\cdot,t)+TV\theta(\cdot,t)
\le C\sigma(t+1)^{-\frac14}+C\delta e^{-\nu t},\quad t\in[0,T).
\end{equation}
Again, taking $\sigma_0$ and $\delta_0$ small we apply Lemma \ref{lemma} to extend the solution $(v,\tilde u)(\cdot,t)$ of \eqref{KS2}, \eqref{ID2} to $t\in[0,\infty)$. Thus, we establish the existence and uniqueness of solutions stated in Theorem \ref{main_thm}. Besides, \eqref{l1uv} and \eqref{4.102} now hold for $0\le t<\infty$. They give \eqref{main_thm_1} and \eqref{main_thm_1a}. We thus finish the proof of Theorem \ref{main_thm}.

%%%%%%%%%%%%%%%%%%%%
\section*{Acknowledgments}%%%%%%
%%%%%%%%%%%%%%%%%%%

 The research of Y. Zeng was partially supported by the National Science Foundation under grant DMS-1908195.  
The research of G. Chen  was partially supported by the National Science Foundation under grant DMS-2008504 and DMS-2306258.

\end{document}